\documentclass[pdflatex
]{sn-jnl}
\usepackage{geometry}
\usepackage{color}
\usepackage[applemac]{inputenc}
\usepackage{graphicx}
\usepackage{amsmath}
\usepackage{amssymb,amsthm}
\usepackage{float}
\renewcommand{\P}{\mathbb P}
\renewcommand{\d}{\delta}

\newtheorem{prop}{Proposition}[section]

\begin{document}

\title[A multi-season epidemic model]{A multi-season epidemic model with random genetic drift and transmissibility}
\author*[1]{\fnm{Tom} \sur{Britton}}\email{tom.britton@math.su.se}

\author[2]{\fnm{Andrea} \sur{Pugliese}}\email{andrea.pugliese@unitn.it}

\affil*[1]{\orgdiv{Department of Mathematics}, \orgname{Stockholm University}, 
\country{Sweden}}

\affil[2]{\orgdiv{Department of Mathematics}, \orgname{University of Trento}, 
\country{Italy}}



\abstract{
 We consider a model for an influenza-like disease, in which, between seasons, the virus makes a random genetic drift $\delta$, (reducing immunity by the factor $\delta$) and obtains a new random transmissibility $\tau$ (closely related to $R_0$). Given the immunity status at the start of season $k$: $\textbf{p}^{(k)}$, describing community distribution of years since last infection,  and their associated immunity levels $\boldsymbol{\iota}^{(k)}$, the outcome of the epidemic season $k$, characterized by the effective reproduction number $R_e^{(k)}$ and the fractions infected in the different immunity groups $\textbf{z}^{(k)}$, is determined by the random pair $(\delta_k, \tau_k)$. It is shown that the immunity status $(\textbf{p}^{(k)}, \boldsymbol{\iota}^{(k)})$, is an ergodic  Markov chain, which converges to a stationary distribution $\bar \pi (\cdot) $. More analytical progress is made for the case where immunity only lasts for one season. We then characterize the stationary distribution of $p_1^{(k)}$, being identical to $z^{(k-1)}$. Further, we also characterize the stationary distribution of $(R_e^{(k)}, z^{(k)})$, and the conditional distribution of $z^{(k)}$ given $R_e^{(k)}$.  The effective reproduction number $R_e^{(k)}$ is closely related to the initial exponential growth rate $\rho^{(k)}$ of the outbreak, a quantity which can be estimated early in the epidemic season. As a consequence, this conditional distribution may be used for predicting the final size of the epidemic based on its initial growth and immunity status.}

\maketitle

\section{Introduction}

Several infectious diseases are globally endemic but circulate around the globe leading to seasonal outbreaks in the winter seasons in the northern and southern hemisphere, seasonal influenza being the prime example. Coming into a new epidemic season, a given community has some natural immunity originating from people being infected in recent years. Vaccinations, typically among the elderly population, reduce the number of severe cases but only marginally affect disease spreading, and are hence neglected in what follows. Between seasons, the circulating virus has often evolved (drifted),  leading to decay of immunity among previously infected, and the transmissibility of the new strain may also differ from the transmissibility of the previous strain (disease severity may also differ,  but this is not considered here).

One of the first to analyse mathematical models for these type of situations was Viggo Andreasen who wrote an influential paper where he analyses a determinstic model for seasonal outbreaks in which the immunity of individuals in the community decays in a specified deterministic way with the number of years since the last infection (Andreasen, 2003). When considering the special case where immunity only lasts one or two years he characterized the possible scenarios. For immunity lasting only one year he showed that there are three parameter regions in which, respectively, the yearly outbreak size converges to a constant, there is an outbreak only every other year, or the outbreak size alternates between a small and a large outbreak every other year. This model has later been extended in several ways. Kucharski et al (2016) consider a virus appearing in  different variants, possibly having (partial) cross immunity to each other, and where immunity of the different strains depend on time since last infection of this strain.  Roberts et al. (2024) extend the model by assuming a heterogeneous community (specified by a continuous distribution) affecting the immunity and infectivity of the individuals. Additionally they study what happens with the long term behaviour of the season outbreaks if the yearly deterministic immunity drift is perturbed randomly by a small amount each year.

In the current paper we take a more genuinely stochastic perspective by treating the yearly genetic drift, and also a transmissibility, $(\delta_k,\tau_k)$ as a bivariate random vector, and where we assume the drift and transmissibility of different years to be independent. One reason for introducing such randomness is to allow for variable outcomes of seasonal outbreaks (as empirically observed), but also because the drift and transmissibility of a new strain depends on mutations in many loci on the virus, their common effect being highly complicated and thus conveniently modelled by random variables. The drift and transmissibility a given year may very well depend on each other, perhaps most likely having negative correlation (one of them being large may be enough to become dominant). Another novel feature in our model is that we allow also the transmissibility to vary between seasons. As a consequence, two quantities of the seasonal outbreaks will vary: the initial growth \emph{and} the final epidemic size, as also observed empirically \cite{C07}. Considering a stochastic model also affect the analysis and in particular the resulting yearly epidemic outcome. Rather than converging to a fixed point or a deterministic periodic situation as in the deterministic setting, the stochastic model leads to a yearly epidemic outbreak depending on previous outbreaks but still being random, having some limiting stationary distribution depending on the random vector $(\delta,\tau)$.

In Section \ref{Sec-Model} we define the seasonal epidemic model with random genetic drift and transmissibility, and also characterize it as a Markov chain with a stationary distribution. In Section \ref{Sec-r=2} we treat the case where immunity only lasts one season, for which more analytical progress is made. In Section \ref{Sec-Numerics} we present simulations for both when immunity lasts one or several years. The paper ends with a short discussion.

\section{The multi-season model with random genetic drift and transmissibility}\label{Sec-Model}

\subsection{Model definition}

The seasonal epidemic is driven by the yearly random genetic drift $\delta_k$ and the yearly random transmissibility $\tau_k$, where the random pairs $(\delta_1, \tau_1), (\delta_2, \tau_2), \dots $ are assumed independent and identically distributed, but where $\delta_k$ and $\tau_k$ may depend on each other (in applications possibly being negatively correlated).

We assume that the immunity of an individual coming into season $k$ depends on how many years $j$ ago the individual was last infected, and that this immunity $\iota_j^{(k)}$ is given by the previous $j-1$ genetic drifts as follows: 
$$
\iota_1^{(k)}=1,\qquad \iota_j^{(k)}=(1-\delta_{k-1})\cdots (1-\delta_{k-j+1}),\ 2\le j<r,\qquad \text{ and } \iota_r^{(k)}=0,
$$
meaning that immunity decays multiplicatively year by year, each year proportionally to the genetic drift. We make an additional assumption: that all immunity is lost after $r\ge 2$ years, so $\iota_r^{(k)}=0$ for all $k$ irrespective the genetic drifts.

What happens during outbreak season $k$ depends on the drift and transmissibility $ (\delta_k, \tau_k)$ but also on the community distribution of when individuals were last infected: we introduce $\textbf{p}^{(k)}=(p_1^{(k)}, \dots , p_{r}^{(k)})$, where $p_j^{(k)}$ denotes the community fraction which were last infected $j$ seasons ago, and $p_r{(k)}=1-\sum_{j=1}^{r-1} p_j^{(k)}$ includes also those who were infected more than $r$ seasons ago, or never: all of them lacking any immmunity.

Given the vector $\textbf{p}^{(k)}$ and their associated immunity levels $\boldsymbol{\iota}^{(k)}=(\iota_{1}^{(k)},\dots ,\iota_{r-1}^{(k)})$ and the drift and transmissibility of year $k$: $(\delta_k,\tau_k)$, the effective reproduction number is given by
\begin{equation}
\label{R_E_k}
R_E^{(k)}=\tau_k \sum_{j=1}^{r} p_j^{(k)} \left(1-(1-\delta_k)\iota_j^{(k)}\right) .
\end{equation}
The basic reproduction number, corresponding to the situation where there is no immunity, is simply $R_0^{(k)}=\tau_k$.
Further, the final epidemic outcome season $k$ in the different immunity classes, $\textbf{z}^{(k)}=(z_1^{(k)},\dots , z_r^{(k)})$, is determined by the drift and transmissibility year $k$: $(\delta_k,\tau_k)$ through the following final size equations

\begin{align*}
1-z_1 & = e^{-(1-(1-\delta_{k})) \tau_k (\sum_j  p_j^{(k)} z_j )  }
\\
\dots &
\\
1-z_j  & = e^{-(1-(1-\delta_{k})\iota_j^{(k)}) \tau_k (\sum_j  p_j^{(k)} z_j )  }
\\
\dots &
\\
1-z_r &= e^{-\tau_k(\sum_j  p_j^{(k)} z_j )} ,
\end{align*}
It is well known that these equations have a unique positive solution $\textbf{z}^{(k)}$ iff $R_E^{(k)} > 1$ (e.g.\ Sec 13.2.3 in \cite{DHB13}). The component $z_j^{(k)}$ denotes the fraction infected during epidemic season $k$ among those who were last infected $j$ seasons ago. The overall fraction getting infected is hence $z^{(k)}=\sum_jz_j^{(k)}p_j^{(k)}$. Summing the equations above, weighted by $p_j^{(k)}$, gives the following equation for the overall fraction getting infected
\begin{equation}
\label{inf_fr}
    1-z =\sum_j p_j^{(k)}e^{-(1-(1-\delta_{k})\iota_j^{(k)}) \tau_k z}.
\end{equation}
This is a 1-dimensional equation for $z^{(k)}$, the overall fraction getting infected, from which the fraction infected of the different types is easily obtained: $z_j^{(k)} = 1- e^{-(1-(1-\delta_{k})\iota_j^{(k)}) \tau_k z^{(k)}}$.

\subsection{Characterization as a Markov chain and its stationary distribution}
We now describe the epidemic model as a Markov chain. Given the integer $r\ge 2$ indicating how many years until immunity is completely lost, the state of our Markov chain is given by the vector $(\textbf{p}^{(k)}, \boldsymbol{\iota}^{(k)})$, where $\textbf{p}^{(k)}=(p_1^{(k)},\dots p_r^{(k)})$ describes the community distribution of time since last infection, and $\boldsymbol{\iota}^{(k)}$ their associated immunity levels $\boldsymbol{\iota}^{(k)}=(\iota_{1}^{(k)},\dots ,\iota_{r-1}^{(k)})$, with $\iota_j^{(k)}=(1-\delta_{k-1})\cdot \dots \cdot (1-\delta_{k-(j-1)})$. The vector   $(\textbf{p}^{(k)}, \boldsymbol{\iota}^{(k)})$ in fact only has dimention $2r-3$, since $p_r^{(k)}=1-\sum_{j=1}^{r-1}p_j^{(k)}$ and $\iota_1^{(k)}=1$.

We now show that the vector $\{ (\textbf{p}^{(k)}, \boldsymbol{\iota}^{(k)})\}_{k=1}^\infty$ is a Markov chain driven by $(\delta_k,\tau_k)$. Given $(\textbf{p}^{(k)}, \boldsymbol{\iota}^{(k)})$, the next state in the Markov chain, $(\textbf{p}^{(k+1)}, \boldsymbol{\iota}^{(k+1)})$, is determined by $(\delta_k,\tau_k)$ as follows. The final outcome of the epidemic season $k$: $\textbf{z}^{(k)}=(z_1^{(k)},\dots , z_r^{(k)})$ is determined as described in the previous subsection (if $R_E^{(k)}\le 1$ then all $z_j^{(k)}=0$ and otherwise they are all strictly positive). The fraction having $j$ years since last infection coming into the next epidemic season $k+1$ is then equal to the fraction having $j-1$ years since last infection in season $k$ and who avoid getting infected in season $k$: $p_j^{(k+1)}=p_{j-1}^{(k)}(1-z_{j-1}^{(k)})$. The first and last components are different:  $p_r^{(k+1)}= p_{r}^{(k)}(1-z_{r}^{(k)}) + p_{r-1}^{(k)}(1-z_{r-1}^{(k)})$, and $p_1^{(k+1)}=\sum_{j=1}^r p_j^{(k)}z_j^{(k)}$ (=all who got infected season $k$). The immunity levels are simply shifted one step: $\iota_j^{(k+1)}=\iota_{j-1}^{(k)}(1-\delta_k)$, and $\iota_1^{(k+1)}=1$ and $\iota_{r}^{(k)}=0$.

We have hence shown that $(\textbf{p}^{(k+1)}, \boldsymbol{\iota}^{(k+1)})$ is determined by $(\delta_k,\tau_k)$ given $(\textbf{p}^{(k)}, \boldsymbol{\iota}^{(k)})$, implying that $\{ (\textbf{p}^{(k)}, \boldsymbol{\iota}^{(k)})\}_{k=1}^\infty$ is a Markov chain. A Markov chain is often described by its transition kernel $\pi(\cdot |\cdot )$. It is in principle possible to determine the transition matrix $\pi( (\textbf{p}^{(k+1)}, \boldsymbol{\iota}^{(k+1)})| (\textbf{p}^{(k)}, \boldsymbol{\iota}^{(k)}))$, at least numerically, although it is quite complicated for general $r$. In Section \ref{Sec-r=2} we make more analytical progress for the case $r=2$. Instead we end this section with the following theorem for our Markov chain, valid for arbitrary $r$. 

\begin{prop}
\label{prop_stat}
    Assume the random vector $(\delta,\tau)$ has a continuous density function with strictly positive support on $[0,1)\times [0,\infty)$, and with an atom at $\delta=1$ (no immunity) with positive probability.  Then the Markov chain $\{ (\textbf{p}^{(k)}, \boldsymbol{\iota}^{(k)})\}_{k=1}^\infty$ is recurrent, aperiodic, and converges to a stationary distribution $\bar \pi (\cdot)$.
\end{prop}
\begin{proof}. We first prove recurrence. We do this by showing that any state $(\textbf{p}, \boldsymbol{\iota})$ can reach the state $(0,\dots , 1, 1, 0,\dots, 0)$ with $p_r=1$ and $\iota_1=1$ and all other components being 0 (so no immunity). This will happen if the next $r$ random pairs $(\delta_j,\tau_j)$ have $\delta_j=1$ and $\tau_j<1$, an event with positive probability ($\tau=R_0<1$ imply no outbreak). 

From state $(0,\dots , 1, 1, 0,\dots, 0)$ we now show that any state $(\textbf{p}, \boldsymbol{\iota})$ can be reached after $r$ seasons. We do this by induction in $r$. Start with $r=2$. We then want to find a $\tau$ such that the final outcome $z=p_1$ for some given value $p_1$, starting with no prior immunity. This is easy: we simply choose $\tau_1$ such that $1-p_1=e^{-p_1\tau_1}$, i.e.\ $\tau_1=-\ln(1-p_1)/p_1$.

For the induction step we want to show that we can reach any fixed state $(p_1,\dots p_r,\iota_1=1, \dots ,\iota_{r-1})$ after $r$ seasons, starting with no immunity season 0. We use that, by induction, we can reach any state $(p'_1,\dots p'_{r-1},\iota'_1=1, \dots ,\iota'_{r-2})$ after $r-1$ seasons. Our task is hence to specify a vector $(p'_1,\dots p'_{r-1},\iota'_1=1, \dots ,\iota'_{r-2})$ and a drift and transmissibility $(\delta,\tau)$ which result in the desired $(p_1,\dots p_r,\iota_1=1, \dots ,\iota_{r-1})$. 

We start with the immunities. Set $\delta=1-\iota_2$, and $\iota'_{j}=\iota_{j+1}/(1-\delta),\ j=1,\dots ,r-2$, leading to the desired immunities. 

We now want to specify the vector $(p'_1,\dots p'_{r-1})$ and transmissibility $\tau$ which, together with immunity vector $(\iota'_1=1, \dots ,\iota'_{r-2})$ result in the vector $(p_1,\dots p_r)$ next year. 

If $p_1=0$, implying that there was no outbreak the previous year, this is easy:  we simply set $p'_{j}=p_{j+1}$, $j=1,\dots ,r-1$, and we choose any $\tau<1$ so that no outbreak will occur. 

The case $p_1>0$ is a bit more involved. 
The vector $(p_1,\dots p_r)$ will be obtained  by iteratively choosing $\tau_i$, $i=1, 2, \dots$ making our vector $(p_1^{(i)},\dots p_r^{(i)})$ approaching $(p_1,\dots p_r)$ closer and closer. 

Start with some arbitrary value larger than $1$, e.g.\ $\tau_1=1.5$. 

The following step is repeated until suitable precision is obtained. For any $i$, set 
$$
1-z'_j=e^{-(1-(1-\delta)\iota_j')\tau_ip_1},\text{ and } p'_{j-1}=p_j/(1-z'_{j-1}), j=1,\dots r.
$$
These candidate values for $\{z'_j\}$ (which depend on $i$ through $\tau_i$!) give the right relations between them (see end of the previous subsection), but it may give too few or too many overall infected. This we now control by Equation (\ref{inf_fr}). Compute $1-z'=\sum_j p'_je^{-(1-\iota'_j) \tau_i p_1}$. If $1-p_1>1-z'$ it means too many got infected overall for our candidate value $\tau_i$, implying that we should pick a smaller $\tau_{i+1}$ and repeat the procedure. On the other hand, if $1-p_1<1-z'$ it means too few got infected overall for our candidate value $\tau_i$, implying that we should pick a larger $\tau_{i+1}$ and repeat the procedure. In each iteration we can pick a new $\tau_i$ being half way between the previous $\tau_{i-1}$ and the closest earlier tested value in the desired direction, implying that our sequence $\tau_i$ will converge, and we stop when our vector $(p'_1,\dots p'_{r-1})$ and $\tau_i$ comes close enough to the desired outcome $(p_1,\dots p_{r})$. This completes the induction step.

The aperiodicity follows from the fact that the immunity-free state $(0,\dots , 1, 1, 0,\dots, 0)$ can return to itself with positive probability (if $\tau<1$).

Since our Markov chain $\{ (\textbf{p}^{(k)}, \boldsymbol{\iota}^{(k)})\}_{k=1}^\infty$ returns to the immunity-free state $(0,\dots , 1, 1, 0,\dots, 0)$ in $r$ steps with positive probability it is a Harris chain (cf.\ Sec 5.8 in Durrett, 2019). It then follows that it has stationary measure $\bar\pi(\cdot)$ to which the chain converges irrespective of starting value. 
\end{proof}

\section{Analytical results for the case $r=2$ where immunity is lost after two seasons}\label{Sec-r=2}
We now focus on the case where $r=2$, so immunity only lasts the year after being infected. In this case our vector $ (\textbf{p}^{(k)}, \boldsymbol{\iota}^{(k)})$ is simply the one dimensional quantity $p_1^{(k)}=:p^{(k)}$, since the outcome on year $k$ only depends on $(\delta_k,\tau_k)$ and how many that were infected last season $p_1^{(k)}=z^{(k-1)}$. We will hence study the distribution of the outcome some year given the outcome of the previous year. We drop the $k$-notation and hence consider the outcome of $R_E$ and $z$, given $p$. All results below are hence conditional on $p$, the fraction infected in the previous season.

We assume that joint probability distribution of the pair $(\delta,\tau)$ has density $q$ (with support in $(0,1) \times (0,+\infty)$), and we wish to explicitly compute some properties of the Markov chain in this specific case.

\subsection{Bivariate distribution of effective reproduction ratio and attack ratio}
Two important properties of a seasonal epidemic are its  effective reproduction ratio $R_e$ and its final attack ratio $z$. The value of $R_e$ determines the initial exponential growth rate, and can be estimated in the first phase of an epidemic. The final attack ratio $z$ (i.e.\ the fraction of the population that gets infected in one season) determines (together with the infection severity) the impact of an epidemic; by definition, it can be ascertained only after an epidemic is over, and may be difficult to estimate, because of under-reporting and other concurring infections. 

As stated above, when $r=2$, the Markov chain is one-dimensional, and only depends on the fraction $p$ infected in the previous year. Simplifying for the case $r=2$ the formulae \eqref{R_E_k} and \eqref{inf_fr}, we have
\begin{equation}
\label{eff_R}
R_e = \tau (p \delta + 1 -p).
\end{equation}
If $R_e \le 1$, then $z = 0$ is the only solution, but if $R_e > 1$, there is a unique solution $z$ in $(0,1)$ to the equation
\begin{equation}
\label{eq_z}
1-z = p e^{-\tau\delta z} + (1-p) e^{-\tau z}.
\end{equation}

In order to consider the probability distribution of the pair $(z,R_e)$, given $p$ (the fraction infected in the previous year), we need to distinguish the case $p=0$ from $p>0$.
We obtain the following
\begin{prop}
    \label{prop_biv}
    Assume that the pair $(\delta,\tau)$ has density $q(\cdot, \cdot)$ (with support in $(0,1) \times (0,+\infty)$.\\
    If $p^{(k)}=0$, the distribution of the pair $(z^{(k)},R_e^{(k)})$ lies on the one-dimensional curve $(z(R_e),R_e)$ where, if $R_e > 1$, $z(R_e)$ is the unique solution in (0,1) of $1-z=e^{-R_e z}$, while $z(R_e)=0$ if $R_e \le 1$. \\ On that curve $R_e^{(k)}$ has density $q_R(y)=\int_0^1 q(x,y)\, dx$.\\
    If $p^{(k)}=p > 0$, the distribution of $(z^{(k)},R_e^{(k)}$ is composed of two parts. With probability
 \begin{equation}
\label{Pz=0}
\P(z=0|p) = \P(R_e \le 1|p) = \int_{ 0}^{ 1}\;\int_{ 0}^{ 1/(p \delta + 1 -p)} q(\delta,\tau)\,d\tau\;d\delta.
\end{equation}
the pair lies in the segment $\{ z=0,\ 0<R_e \le 1\}$  with density $$\int_0^1 q\left(\delta,\frac{R_e}{p\delta+1-p}\right)\frac{1}{p\delta+1-p}\,d\delta.$$
The second part lies in $\{ 0<z<1,\ R_e >1\}$ with bivariate density
 \begin{equation}
\label{biv_dens2}
f(z,R_e|p) = \begin{cases} q(\delta^*,\tau^*) \frac{ p \delta^* + 1 -p-R_e\left(p\delta^* e^{- \delta^* \tau^* z } + (1-p) e^{-  \tau^* z  }\right) }{ p(1-p)R_e z\left(e^{- \delta^* \tau^* z }-e^{-  \tau^* z }\right)}&\mbox{if }(z,R_e) \in \mathcal{A}\\
0 &\mbox{otherwise } \end{cases}
\end{equation}
where
\begin{multline}
\label{range_F}
\mathcal{A}= \{0<z<1-p,\ - \frac{\log(1-z)}{ z} < R_e < - \frac{(1-p)}{ z} \log\left(1 - \frac{z}{1-p }\right)\}\\
\cup
\{1-p \le z < 1,\  - \frac{\log(1-z)}{ z} < R_e\},
\end{multline}
$\delta^*(z,R_e)$ is the unique solution (existing when $(z,R_e) \in \mathcal{A}$) in $(0,1)$ of 
\begin{equation}
\label{function_G}
G(\delta,z,R_e) := p e^{- \frac{\delta R_e z}{ p \delta + 1 -p} } + (1-p) e^{- \frac{ R_e z}{ p \delta + 1 -p} } + z - 1 = 0.
\end{equation}
and
$$ \tau^* = \frac{R_e}{ p \delta^* + 1 -p}.$$
\end{prop}
\begin{proof}
    If $p=0$, formulae \eqref{eff_R} and \eqref{eq_z} show that $R_e=R_0=\tau$ and $z$ is the only solution of $1-z=e^{-\tau z}$, which is completely determined by $\tau$. Hence the pair $(z,R_e)$ lies on a one-dimensional curve: $R_e$ has density $q_R(y)=\int_0^1 q(x,y)\, dx$, while $z$ is a function of $R_e$. The density of $R_e$ is the density of $\tau$ which is obtained by marginalizing the density $q$ as in the statement of the Theorem.

On the other hand, if $p >0$, we can consider the map
\[\begin{split} F : (0,1) \times (0,+\infty) &\to [0,1) \times (0,+\infty)\\
 (\delta,\tau)\ &\to\ (z,R_e). \end{split}\]
 Using the standard formula (see. for instance, \cite[Theorem 12.7]{JacodProtter}) for the density of a random variable $F(X)$, the density of the pair $(z,R_e)$ is
 \begin{equation}
\label{biv_dens1}
f(z,R_e) = \begin{cases} q(F^{-1}(z,R_e))\times  |\mbox{det}((F^{-1})'(z,R_e))|&\mbox{if }(z,R_e) \in \mbox{Im}(F)\\
0 &\mbox{if }(z,R_e) \not\in \mbox{Im}(F). \end{cases}
\end{equation}
In this case, $F_1((\delta,\tau))\equiv 0$ for all values such that $F_2((\delta,\tau))\le 1$; thus $F$ is not invertible everywhere. It is useful to think of \eqref{biv_dens1} as a defective density restricted to the set $\{(z>0,\ R_e > 1)\}$, while, when $R_e \le 1$ the support of $(z^{(k+1)},R_e^{(k+1)})$ lies in the segment with $z=0$. It is easy to see that $\P(z=0) = \P(R_e \le 1)$ can be computed as in \eqref{Pz=0}, while on that segment the density of $R_e^{(k+1)}$ can be computed from 
$$
\P(R_e^{(k+1)}\le R_e) = \int_{ 0}^{ 1}\;\int_{ 0}^{ R_e/(p \delta + 1 -p)} q(\delta,\tau)\,d\tau\;d\delta.
$$

Considering now \eqref{biv_dens1},
the density $f(z,R_e)$ cannot be written explicitly, since there are not explicit expressions for $F$ and $F^{-1}$. However, one can improve on \eqref{biv_dens1} by computing the range of $F$ and the derivatives of $F^{-1}$.

In order to find the range of $F$, we first find
$$ \tau = \frac{R_e}{ p \delta + 1 -p}$$
and then, substituting in \eqref{eq_z}, $\delta$ can be obtained as a solution of $G(\delta,z,R_e)=0$ with $G$ shown in \eqref{function_G}.

Simple computations then yield
\begin{align}
\label{derG1} \frac{\partial G}{\partial \delta } &= - \frac{p(1-p)R_e z}{( p \delta + 1 -p)^2 }\left(e^{- \frac{\delta R_e z}{ p \delta + 1 -p} }-e^{- \frac{ R_e z}{ p \delta + 1 -p} }\right) < 0 \\
\label{derG2}\frac{\partial G}{\partial R_e } &= - \frac{z}{ p \delta + 1 -p}\left(p\delta e^{- \frac{\delta R_e z}{ p \delta + 1 -p} } + (1-p) e^{- \frac{ R_e z}{ p \delta + 1 -p} }\right) < 0 \\
\label{derG3}\frac{\partial G}{\partial z } &= 1 - \frac{R_e}{ p \delta + 1 -p}\left(p\delta e^{- \frac{\delta R_e z}{ p \delta + 1 -p} } + (1-p) e^{- \frac{ R_e z}{ p \delta + 1 -p} }\right) >0.
\end{align}
The sign of $\frac{\partial G}{\partial z }$ comes from the fact that the derivative of the RHS of \eqref{eq_z} is larger than $-1$ at the solution $z \in (0,1)$ of \eqref{eq_z}.

Since $G$ is decreasing a function of $\delta$, there exists a solution $\delta^*(z,R_e)$ of $G(\delta,z,R_e)=0$ if and only if
\begin{equation}
\label{cond_inv}
p+ (1-p)e^{- \frac{ R_e z}{  1 -p} }+z-1  = G(0,z,R_e)  > 0 > G(1,z,R_e) = e^{-  R_e z}+z-1.
\end{equation}
It is easy to see that, if $ z < 1-p $, \eqref{cond_inv} is equivalent to
\[
- \frac{\log(1-z)}{ z} < R_e < - \frac{(1-p)}{ z} \log\left(1 - \frac{z}{1-p }\right), \qquad\]
while only the first condition holds if $z \ge 1-p$.

It follows that
\[
\mbox{Im}(F) = \mathcal{A}
\cup\{z=0,\ R_e\le 1\},
\]
with $\mathcal{A}$ specified in \eqref{range_F}.

Considering only the case where $R_e > 1$, where $F$ is invertible and \eqref{biv_dens1} holds, we have
$$ F^{-1}(z,R_e) = \left(\delta^*(z,R_e), \frac{R_e}{ p \delta^*(z,R_e) + 1 -p}\right). $$
We can then compute its Jacobian using \eqref{derG1}--\eqref{derG3}. First of all, notice that
$$ \frac{\partial \tau}{\partial z } = - \frac{pR_e}{(p \delta^* + 1 -p)^2 } \frac{\partial \delta}{\partial z } \qquad \frac{\partial \tau}{\partial R_e } = \frac{1}{ p \delta^* + 1 -p}- \frac{pR_e}{(p \delta^* + 1 -p)^2 } \frac{\partial \delta}{\partial R_e }. $$
Then
$$
\left| \begin{array}{cc}\dfrac{\partial \delta}{\partial z }& \dfrac{\partial \delta}{\partial R_e } \\[1em]
\dfrac{\partial \tau}{\partial z } & \dfrac{\partial \tau}{\partial R_e }\end{array} \right| = \frac{\partial \delta}{\partial z } \cdot \frac{1}{ p \delta^* + 1 -p}. $$
Using the implicit function theorem, and remembering $ \tau = \frac{R_e}{ p \delta + 1 -p}$,
$$ \frac{\partial \delta}{\partial z } = - \frac{ \frac{\partial G}{\partial z } }{  \frac{\partial G}{\partial \delta } }= ( p \delta^* + 1 -p)\frac{ p \delta^* + 1 -p-R_e\left(p\delta^* e^{- \delta^* \tau^* z } + (1-p) e^{-  \tau^* z  }\right) }{ p(1-p)R_e z\left(e^{- \delta^* \tau^* z }-e^{-  \tau^* z }\right)}. $$
Substituting, we obtain \eqref{biv_dens2}. 
\end{proof}

\subsection{Transition probabilities}
From \eqref{biv_dens2} one can compute several other formulae of interest. We start with the transition probabilities. Under the assumption that $(\delta,\tau)$ has density $q$ (with support in $(0,1) \times (0,+\infty))$, one sees that $p^{(k)}$ has a distribution with an atom at 0 and a density over $(0,1)$, and the same will be true for $p^{(k+1)}.$

Then one obtains the following
\begin{prop}
\label{prop_trans}
    The transition probability of the Markov chain consist of
    \begin{equation}
\label{prob_trans00}
  \P(p^{(k+1)} =0 | p^{(k)} =0) =\pi_{00}= \P(\tau \le 1) = \int_0^1\; \int_0^1 q(x,y)\, dx\, dy. 
\end{equation}
\begin{multline}
\label{prob_trans10}
 \P(p^{(k+1)}=0 | p^{(k)} =p)  =\pi(p,0) =\P(R_e^{(k)} \le 1) = \P(\tau(\d p+(1-p)) \le 1) \\= \int\limits_0^1\; \int\limits_0^{1/(xp+(1-p))} q(x,y)\, dy\, dx. 
\end{multline}
while the density function of $p^{(k+1)}$ for $p'>0$ conditional on $p^{(k)} =p \ge 0$ is
\begin{equation}
\label{prob_trans2}
\pi(p,p') = \int\limits_{ 0}^{1 } q(x,\tau^*)   \frac{1 - \tau^*\left(p x e^{- \tau^* x p'}+(1-p) e^{-\tau^* p'}\right)}{p'\left(px e^{- \tau^* x p'}+(1-p) e^{- \tau^* p'}\right) }\, dx
\end{equation}
where $\tau^*=\tau^*(x,p')$ is the unique solution of
$$ p e^{-x\tau p'} + (1-p) e^{-\tau p'}+ p' - 1 =0. $$
\end{prop}
\begin{proof}
   The computation of \eqref{prob_trans00} and \eqref{prob_trans10} is straightforward.

To compute the density function $\pi(p,p')$ where $p' \in(0,1)$ and $p \ge 0$, one needs only integrate \eqref{biv_dens2} to obtain
\begin{equation}
\label{prob_trans1}
\pi(p,p') = \begin{cases} \int\limits_{ - \frac{\log(1-p')}{ p'} }^{ - \frac{(1-p)}{ p'} \log\left(1 - \frac{p'}{1-p }\right) } f(p',y)\,dy \quad\ & p' < 1-p \\[0.3em]
\int\limits_{ - \frac{\log(1-p')}{ p'} }^{ +\infty } f(p',y)\,dy & p' \ge 1-p.
\end{cases}
\end{equation}
It is more convenient to change the integration variable in \eqref{prob_trans1} to $x = \delta^*(p',y)$; then 
$$dy=  - \frac{ \frac{\partial G}{\partial \delta } }{  \frac{\partial G}{\partial R_e } }dx $$
while
\begin{equation}
\label{f_expl}
f(p',y) = - q(\delta^*(p',y),\tau^*(p',y))  \frac{ \frac{\partial G}{\partial z } }{  \frac{\partial G}{\partial \delta } }\cdot\frac{1}{p \delta^* + 1 -p } . 
\end{equation}

Hence, using \eqref{derG1}--\eqref{derG3},
\begin{multline*}
    f(p',y)\,dy =  - q(x,\tau^*)  \frac{ \frac{\partial G}{\partial z } }{  \frac{\partial G}{\partial R_e } }\cdot \frac{1}{p x + 1 -p }\,dx \\ =  q(x,\tau^*) \frac{1 - \tau^*\left(p x e^{- \tau^* x p'}+(1-p) e^{-\tau^* p'}\right)}{p'\left(px e^{- \tau^* x p'}+(1-p) e^{- \tau^* p'}\right) }\, dx
\end{multline*} 
Substituting in \eqref{prob_trans1}, we obtain \eqref{prob_trans2}.
\end{proof}

\subsection{The conditional distribution of final size given the effective reproduction ratio}
\label{subs:cond_dist}
As mentioned above, coming into season $k$ the fraction infected in the previous season $z^{(k-1)}=p^{(k)}$ is known, and the value of $R^{(k)}_e$ can be estimated from data during the early phase of the epidemic, while the final size $z^{(k)}$ ($=p^{(k+1)}$) can be estimated from data only at the end of a seasonal epidemic. Hence it may be useful to predict $z^{(k)}$ knowing $z^{(k-1)}=p^{(k)}$ and $R^{(k)}_e$.

If $R^{(k)}_e \le 1$, necessarily $z^{(k)}=0$. Similarly if $z^{(k-1)}=0$, there is no immunity; then both $R^{(k)}_e$ and $z^{(k)}$ are determined by $\tau_k$ only, so that knowing one, the other one can be exactly computed.

Finally, the non-trivial case is  where $z^{(k-1)}>0$ and $R_e^{(k}=R_e>1$. We wish to compute the conditional density
\begin{equation}
\label{cond_dens}
f(z^{(k)}=z | R_e^{(k)}=R_e,\ p^{(k-1)}=p ) = \frac{f(z,R_e|p)}{\int\limits_{z': (z',R_e) \in\mathcal{A} }f(z',R_e|p)\, dz'},
\end{equation}
where $f$ denotes the density for whatever variables and conditional events that are indicated. 

In this case too, it is convenient to change the integration variable to $x = \delta^*(p',R_e)$ so that $dp' =  - \dfrac{ \frac{\partial G}{\partial \delta } }{  \frac{\partial G}{\partial z } }dx$. Remembering that $\delta^*$ is the solution of $G(\delta,z,R_e)=0$ as defined in \eqref{function_G}, this yields, together with \eqref{f_expl},
\begin{equation}
\label{cond_dens2}
f(z | R_e,p) = \frac{f(z,R_e|p)}{ \int_0^1 q(x,\tau) \frac{1}{p x + 1 -p }\, dx}
\end{equation}
where in the integral $\tau = \dfrac{R_e}{ p x + 1 -p}$.

\subsection{Stationary distribution}
\label{subs:stat}
The results of Section \ref{Sec-Model} can be made more specific when $r=2$. Precisely
\begin{prop}
Assume that the pair $(\delta,\tau)$ has a density $q$ with support in $(0,1)\times(\tau_m,\tau_M)$ with $0\le \tau_m < 1 < \tau_M \le +\infty$, such that the marginal density $q_\tau(y) = \int_0^1 q(x,y)\,dx > 0$ on $(\tau_m,\tau_M)$. 

Then the Markov chain $\{ p^{(k)}\}_{k=1}^\infty$ is recurrent, aperiodic, and converges to a stationary distribution $\bar \pi $ which consists of an atom $\tilde \pi_0$ at 0 and a density $\tilde \pi(\cdot)$ with support on $(0,\bar z)$; $\bar z = 1$ if $\tau_M = +\infty$, while $\bar z$ is the solution of $1-z=e^{-\tau_M z}$ if $\tau_M < + \infty$.

The  pair $(\bar \pi_0, \bar \pi(x))$  satisfies the system
\begin{equation}
\label{inv_prob}
\begin{split} \bar \pi_0 &= \bar \pi_0 \pi_{00} + \int_0^{\bar z} \pi(x,0) \bar \pi(x)\, dx \\
\bar \pi(x) &= \bar \pi_0 \pi(0,x) + \int_0^{\bar z} \pi(z,x) \bar \pi(z)\, dz\qquad x \in (0,\bar z) \\
1 &= \bar \pi_0 + \int_0^{\bar z}  \bar \pi(x)\, dx.
\end{split}
\end{equation}
\end{prop}
\begin{proof}
    The proof follows the same path as in Proposition \ref{prop_stat}, but much more simply since the state space is one-dimensional. From any state $p$ one reaches $0$ if $\tau < 1$; vice versa from 0 one can reach any state $p \in (0,\bar z)$ for an appropriate $\tau \in (1,\tau_M)$. The aperiodicity follows from the fact that $\pi_{00}= \P(p^{(k+1)} =0 | p^{(k)} =0) >0$. 

It follows that the Markov chain is recurrent, aperiodic, and convergent. Equations \eqref{inv_prob} come directly from the stationarity condition.  
\end{proof}
Presumably, the result holds also if $\tau_m \ge 1$, but the proof would be more delicate. If $\tau_M\le 1$, then $\bar\pi_0=1$, since $P(R_0\le 1)=1$.

\medskip
Similarly to Subsection \ref{subs:cond_dist}, one may wish to predict $z^{(k)}$ knowing $R_e^{(k)}$; in this case we assume that $p^{(k)}$ is not known, but its distribution is the stationary one $\bar \pi$.

One starts by computing the stationary bivariate distribution of $(z^{(k)}, R_e^{(k)})$. If $p^{(k)}=0$ (with probability $\bar \pi_0$), $R_e^{(k)}=\tau_k$ and $z^{(k)} = z(R_e^{(k)})$ where the function $z(\cdot)$ is defined in Proposition \ref{prop_biv}.
Hence, the pair $(z^{(k)}, R_e^{(k)})$ has a component on the curve $(z(R_e),R_e)$ with (one-dimensional) density of $R_e$ given by $\bar\pi_0 q_\tau(\cdot)$, with $q_\tau$ the marginal distribution of $\tau$.

If $p^{(k)}>0$, its (defective) density is given by $\bar \pi(\cdot)$ and one obtains the bivariate distribution simply by multiplying by $\bar \pi(x)$ the formulae presented in Proposition \ref{prop_biv}.

Precisely, the (defective) bivariate density over $\{ z \in (0,1), R_e > 1\}$ is
\begin{equation}
    \label{stat_biv_dens}
    \bar f(z,R_e) = \begin{cases}
        0 & \mbox{if }-\frac{\ln(1-z)}{z}\ge R_e \\
        \int_{1-p(z,R_e)}^1 \bar \pi(x) f(z,R_e | x)\,dx\quad\ & \mbox{otherwise}
    \end{cases}
\end{equation}
where $f(z,R_e | x)$ is given by \eqref{biv_dens2}, and $p(z,R_e)$ is the only solution in $(0,1)$ of 
$$ R_e = - \frac{(1-p)}{ z} \ln\left(1 - \frac{z}{1-p }\right). $$
From the previous formulae, one obtains, through straightforward computations, the conditional stationary distribution of $z^{(k)}$ given $R_e^{(k)}$. One may note that this distribution has an atom at $z(R_e)$ of probability 
$$ \frac{\bar \pi_0 q_\tau(R_e)}{\bar \pi_0 q_\tau(R_e) + \int_0^{z(R_e)} \bar f(z,R_e) dz}$$
where $\pi(x,z)$ is given in \eqref{prob_trans2}, and a continuous (defective) density over $(0,z(R_e))$
$$\bar f(z |R_e) =  \frac{\bar f(z,R_e)}{\bar \pi_0 q_\tau(R_e) + \int_0^{z(R_e)} \bar f(z,R_e) dz}$$.
\section{Numerical illustrations}\label{Sec-Numerics}
We start by illustrating the analytical results obtained in the case $r=2$. We chose 4 different densities $q$ for the pairs $(\delta_k,\tau_k)$; in all cases $\delta$ follows the Beta$(a,b)$ distribution, and $\tau$ the log-normal$(\mu,\sigma^2)$. In cases 1 and 2 the two variables are independent; in cases 3 and 4 they are negatively correlated. Negative correlation between $\delta$ and $\tau$ appears realistic, since viral strains may be fitter either because they are more transmissible or because they are able to evade existing immune response.

\begin{table}[h]
\begin{tabular}{|c|l|l|l|l|}
\hline 
Parameters & Case 1 & Case 2 & case 3 & Case 4\\
\hline 
$a$ & 3 & 4 & 0.5 & 3 \\ 
 $b$ & 7 & 6 & 1.5 & 7 \\ 
 $\mu$ & 0.683  & 1.08  & 0.6 $-$ 0.4 $\delta$ & 0.7 $-$ 0.5 $\delta$ \\ 
 $\sigma^2$ & 0.02 & 0.02 & 0.02 & 0.02 \\ 
 $E(\delta)$ & 0.3 & 0.4 & 0.25 & 0.3 \\ 
 $\sqrt{V(\delta)}$ & 0.14 & 0.15 & 0.25 & 0.14 \\ 
 $E(\tau)$ & 2 & 2.97 & 1.74 & 1.83 \\ 
 $\sqrt{V(\tau)}$ & 0.28 & 0.42 & 0.27 & 0.28 \\ 
 correl.\ coeff.\ & 0 & 0 & -0.38 & -0.32 \\ 
 \hline 
 
\end{tabular}
\caption{Parameters and features of the 4 distributions used in the simulations. In all cases $\delta \sim \mbox{Beta}(a,b)$, while $\tau \sim \mbox{Lognorm}(\mu,\sigma^2)$.}
\label{tab:pairs}
\end{table}

\begin{figure}[h]
    \centering
    \includegraphics[width=0.35 \linewidth]{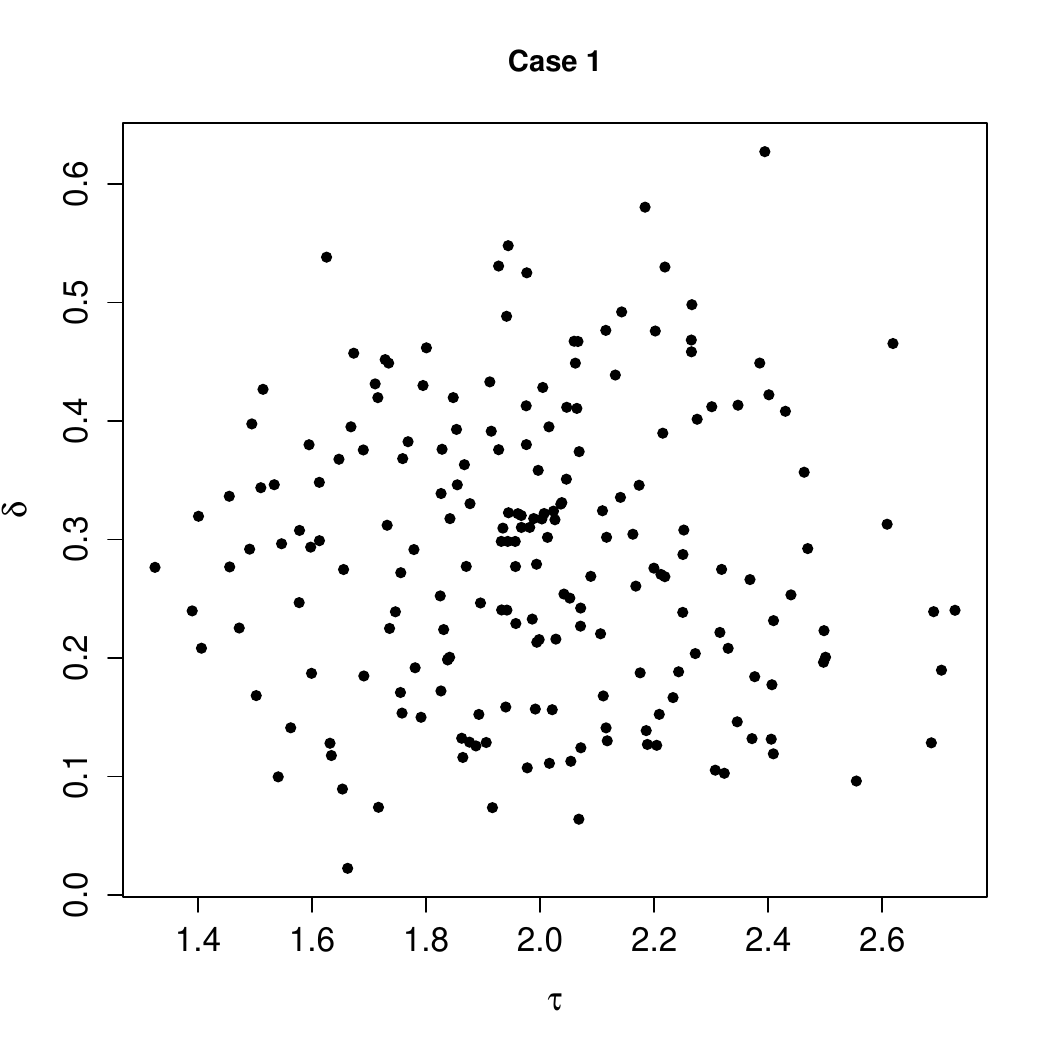}\ 
     \includegraphics[width=0.35 \linewidth]{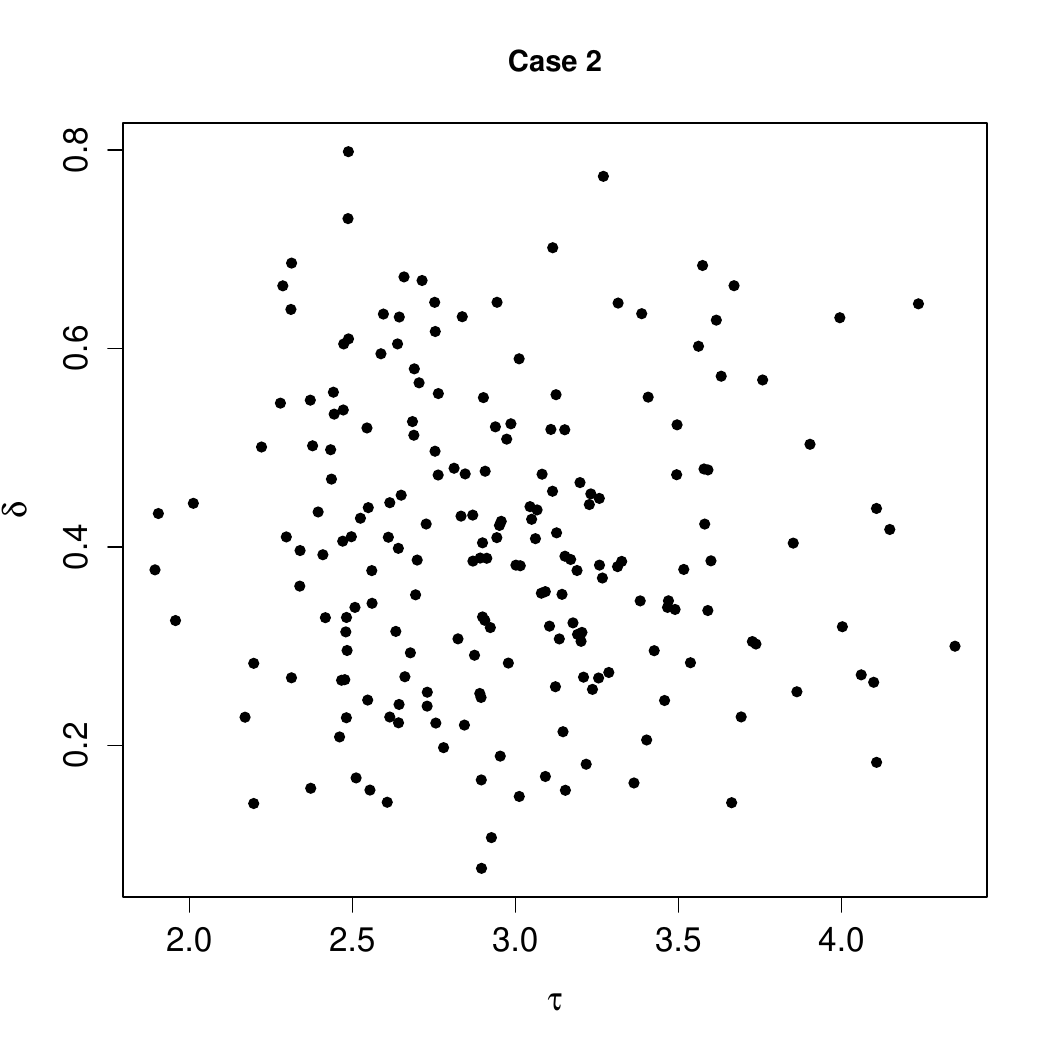}\\
       \includegraphics[width=0.35 \linewidth]{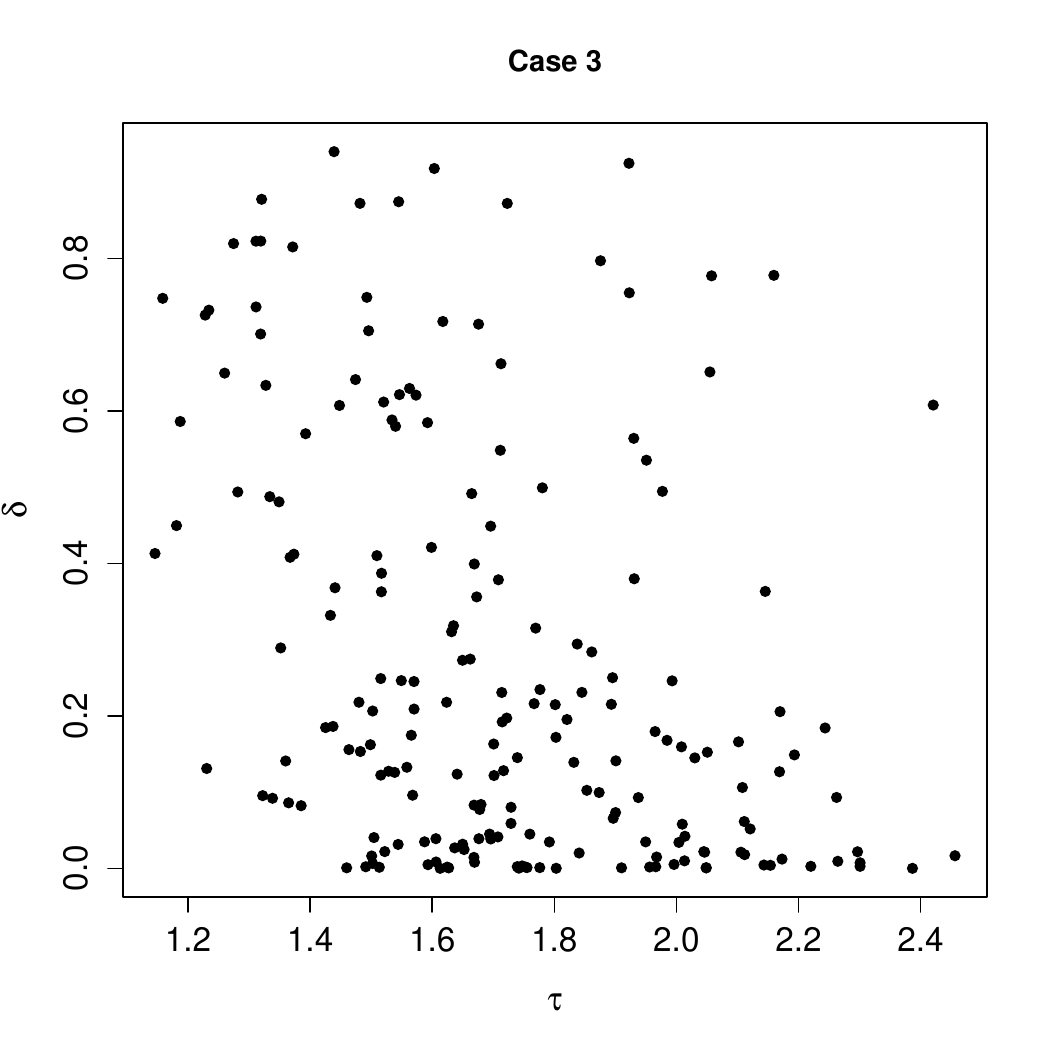}\ 
     \includegraphics[width=0.35 \linewidth]{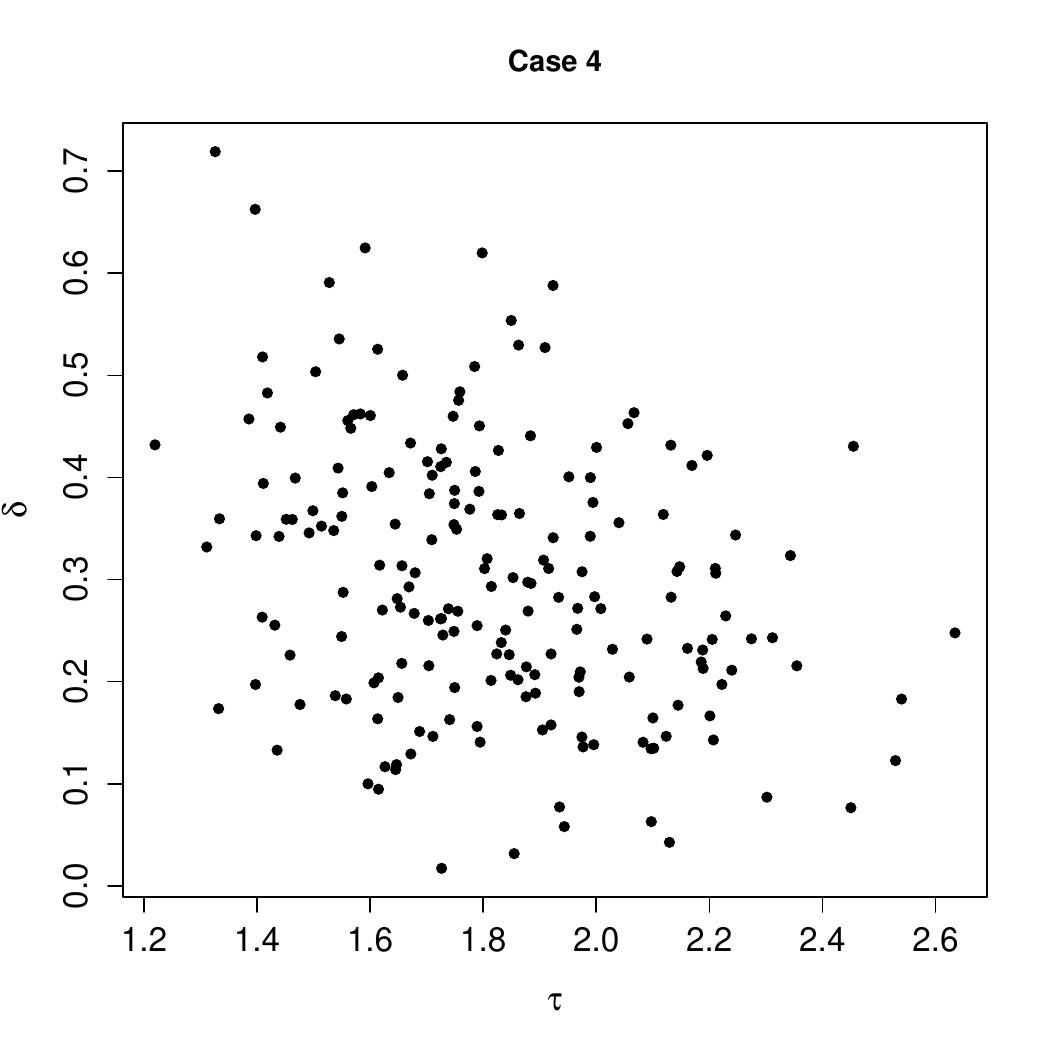}
    \caption{200 draws of the random variables $(\tau,\delta)$ for each case described in Table \ref{tab:pairs}}
    \label{fig_corr}
\end{figure}

Parameters and features of the four cases are listed in Table \ref{tab:pairs}, while Figure \ref{fig_corr} shows 200 random values from each of the four distributions.

\begin{figure}
    \centering
    \includegraphics[width=0.35\linewidth]{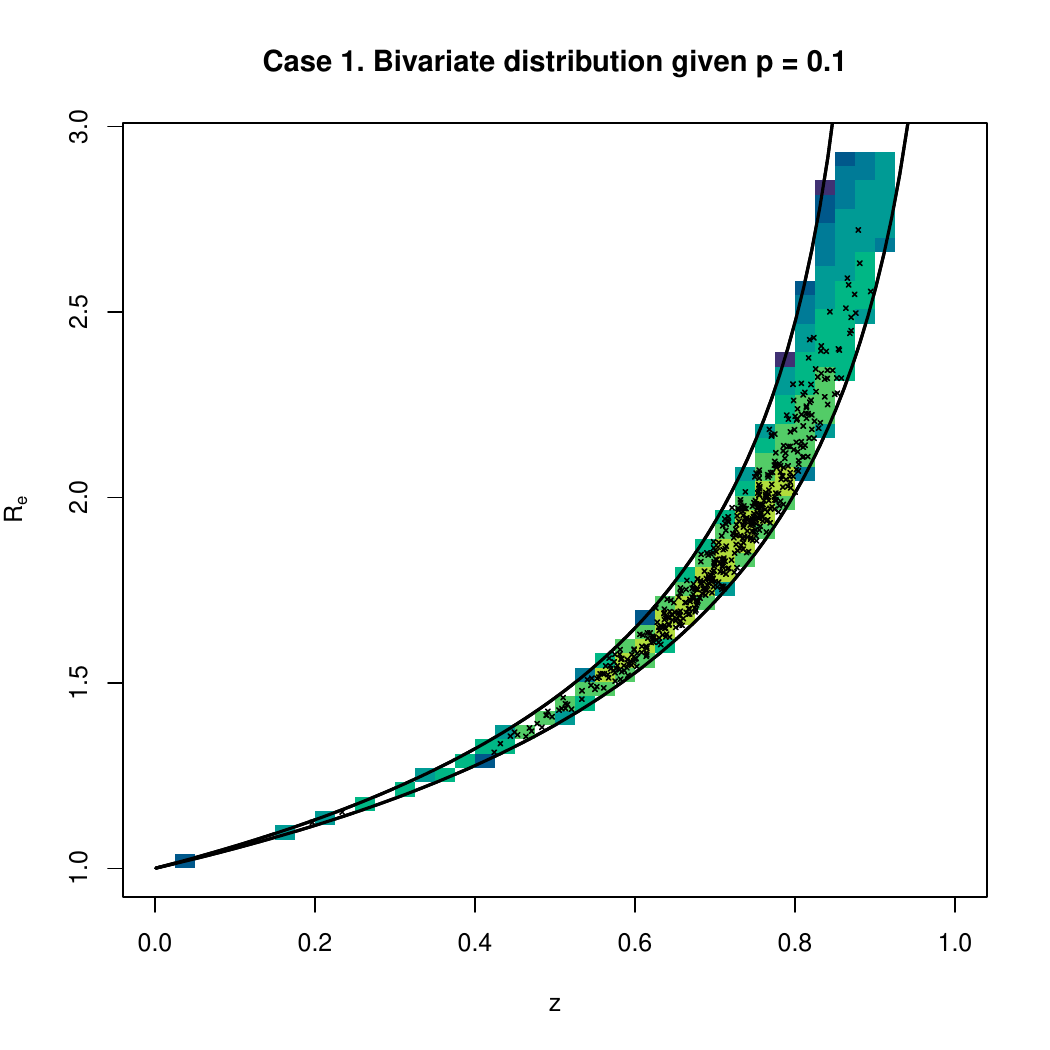}\ 
  \includegraphics[width=0.35\linewidth]{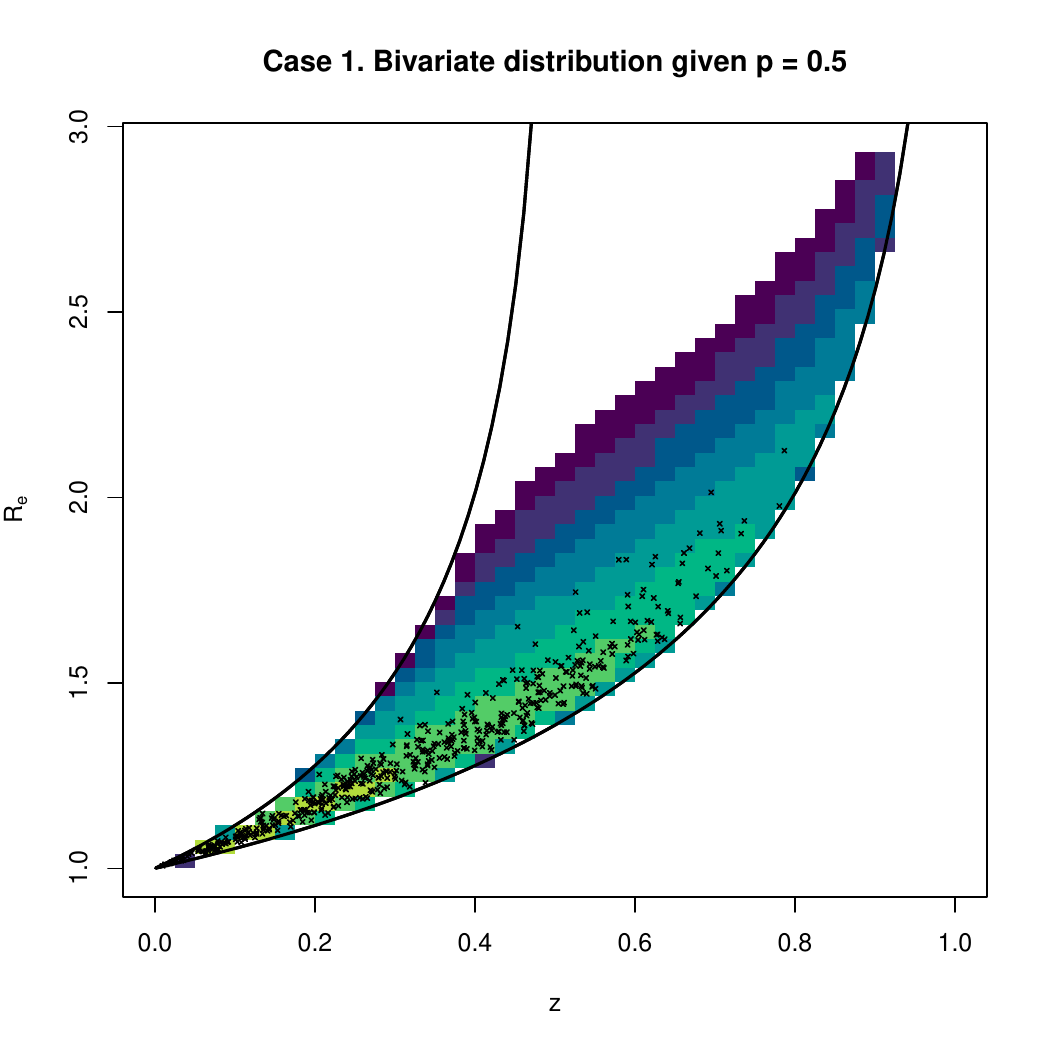}
    \caption{Examples of the bivariate distributions of  $(z^{(k)},R_e^{(k})$ given $p^k = z^{(k-1)} = p$ computed using Proposition \ref{prop_biv} for $p=0.1$ (left panel) or $p=0.5$ (right panel) in case 1 of Table \ref{tab:pairs}.}
    \label{fig:biv_dist}
\end{figure}

In Figure \ref{fig:biv_dist}, we show some heat maps of the bivariate distributions of the pair $(z^{(k)},R_e^{(k)})$ given $p^{(k)} = z^{(k-1)} = p$, computed using Proposition \ref{prop_biv}. To give a better feeling of the distributions we added (small crosses) 500 random draws. It can be seen that either ($p=0.1$), the support of the distribution is very narrow, or, when the support is wider ($p=0.5$), the density is concentrated near the curve $R_e = -\log(1-z)/z$.

\begin{figure}[ht]
    \centering
    \includegraphics[width=0.35 \linewidth]{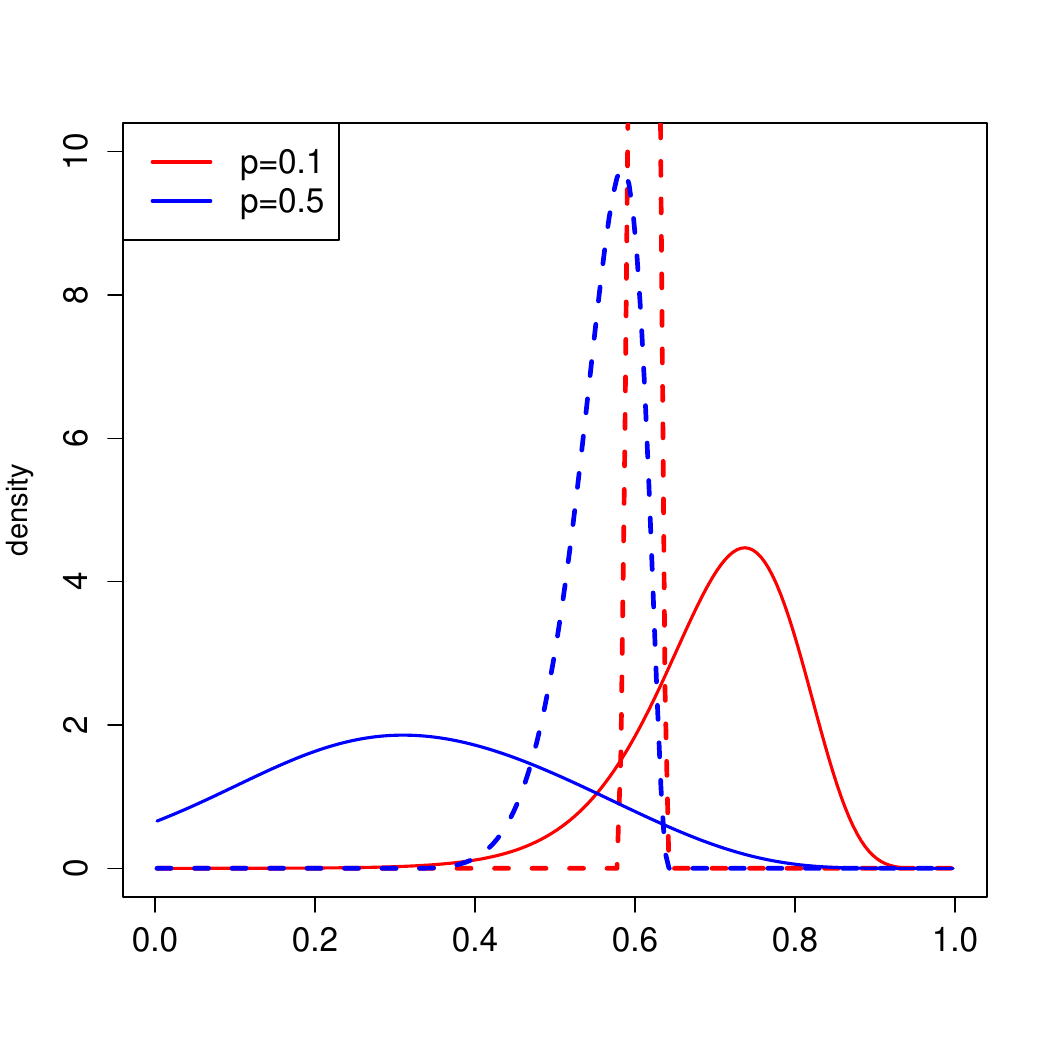}\ 
     \includegraphics[width=0.35 \linewidth]{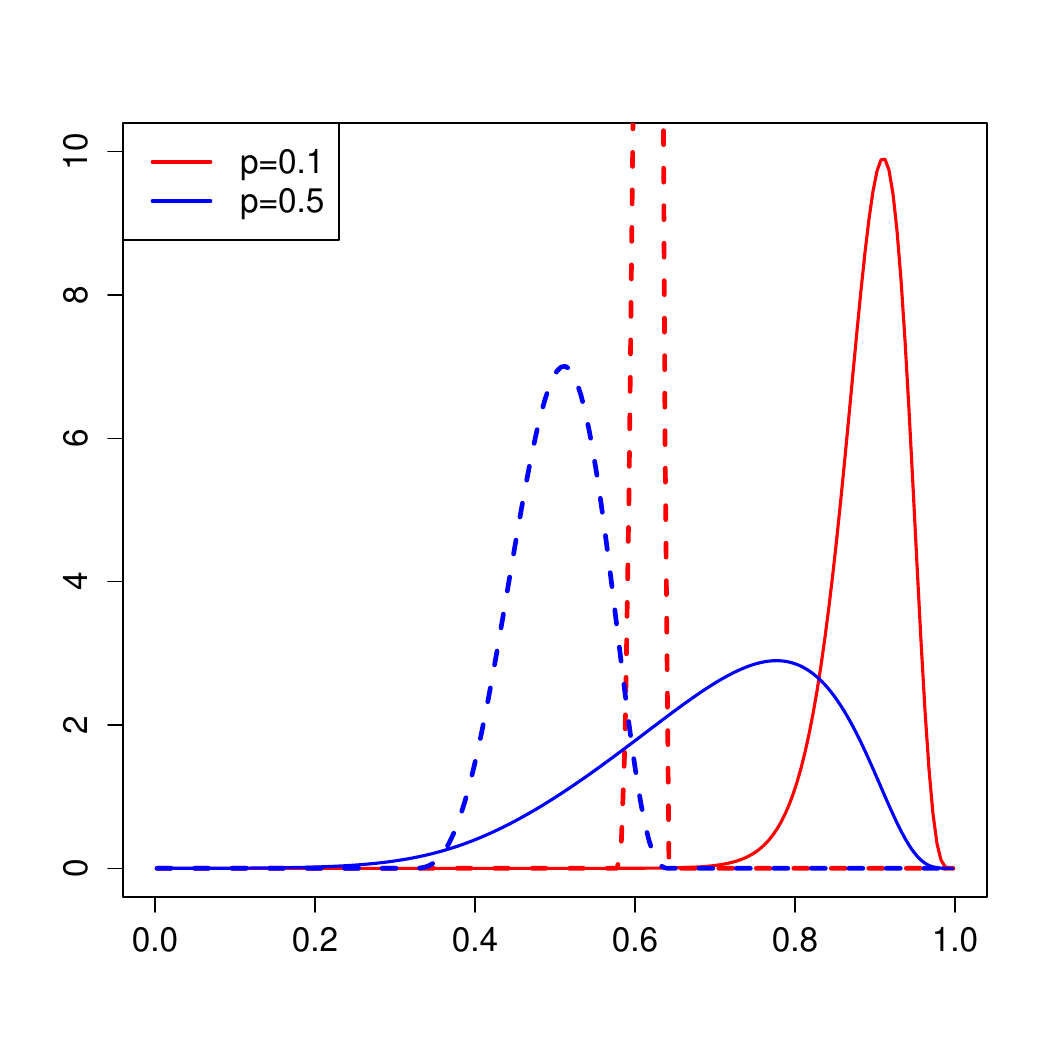}\\[-0.5cm]
       \includegraphics[width=0.35 \linewidth]{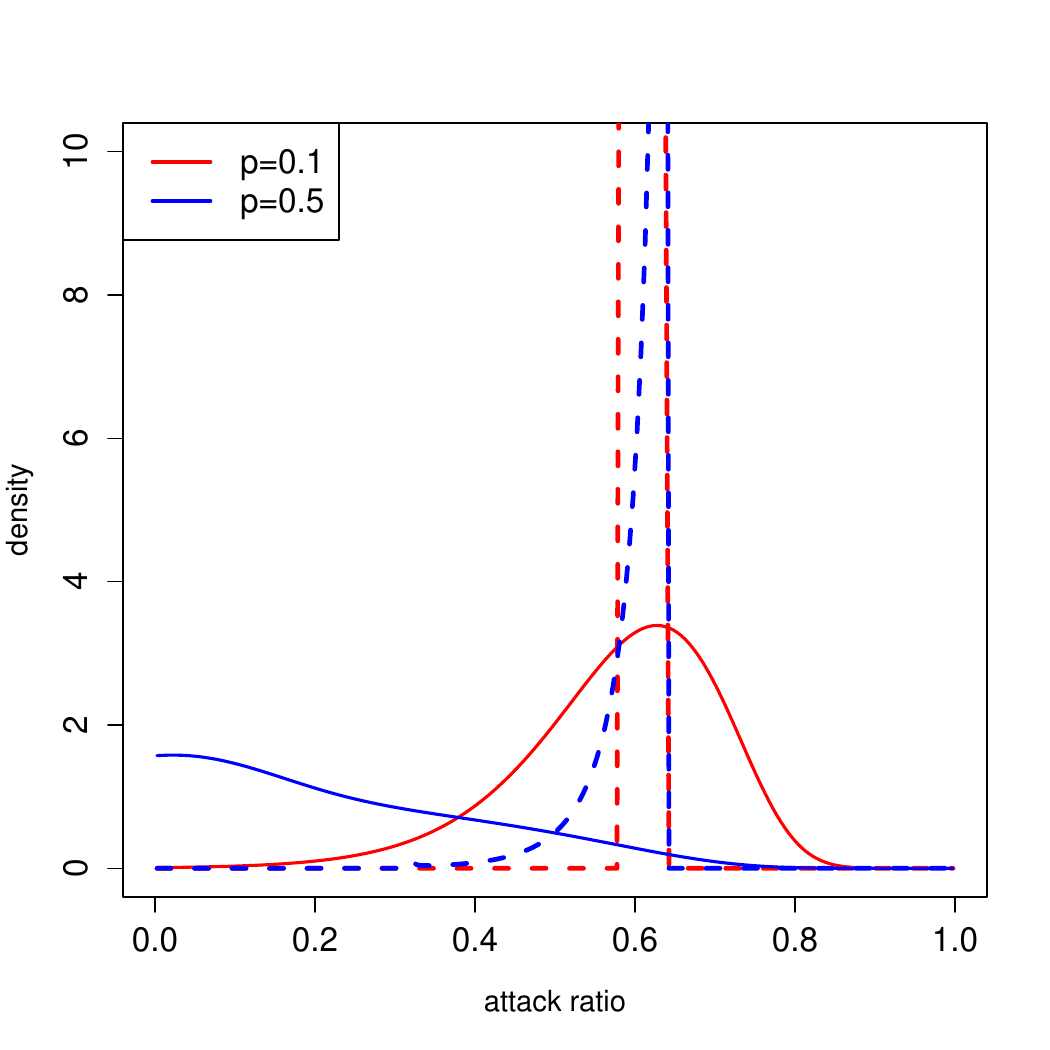}\ 
     \includegraphics[width=0.35 \linewidth]{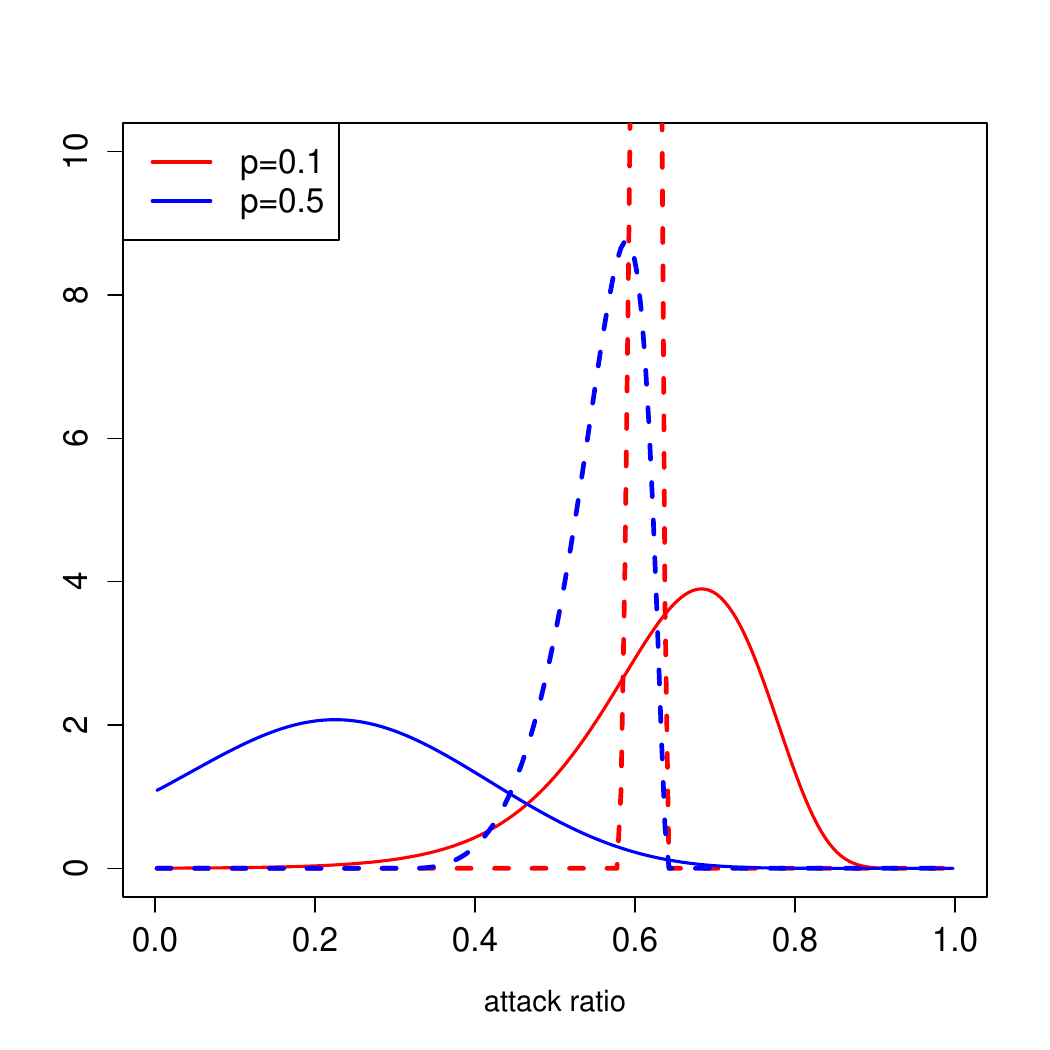}
    \caption{The distributions of $z^{(k)}$ given $p^{(k)}=z^{(k-1)}$  (solid curves), and given the pair $(R_e^{(k)},p^{(k)})$ (dashed curves). The red curves are computed with $p^{(k)}=p=0.1$, the blue ones with $p=0.5$. The conditional distributions are based on $R_e^{(k)}=1.6$ The four panels correspond to the four distributions of the pairs $(\delta_k,\tau_k)$ as shown in Table \ref{tab:pairs}.}
    \label{fig:dist_givenp}
\end{figure}

In Figure \ref{fig:dist_givenp} we show some examples of the distribution of $z^{(k)}= p^{(k+1)}$ given $p^{(k)} = p$ for two different values of $p$; clearly, these correspond to the transition probabilities,  obtained from Proposition \ref{prop_trans}. In the same pictures we also show the distribution assuming that also $R_e^{(k)}$ is known and equal to 1.6. It can be seen that the distributions of $z^{(k)}$ is quite different depending on the value of $p$. Knowledge of 
$R_e^{(k)}$ hence greatly reduces uncertainty (i.e.\ randomness) of $z^{(k)}$. As a consequence, a prediction of the final attack ratio after observing the initial growth in an epidemic will greatly improve the precision of a prediction of the final attack ratio $z^{(k)}$ as compared to only knowing the attack ratio of the year before.

\begin{figure}[h]
    \centering
    \includegraphics[width=0.35 \linewidth]{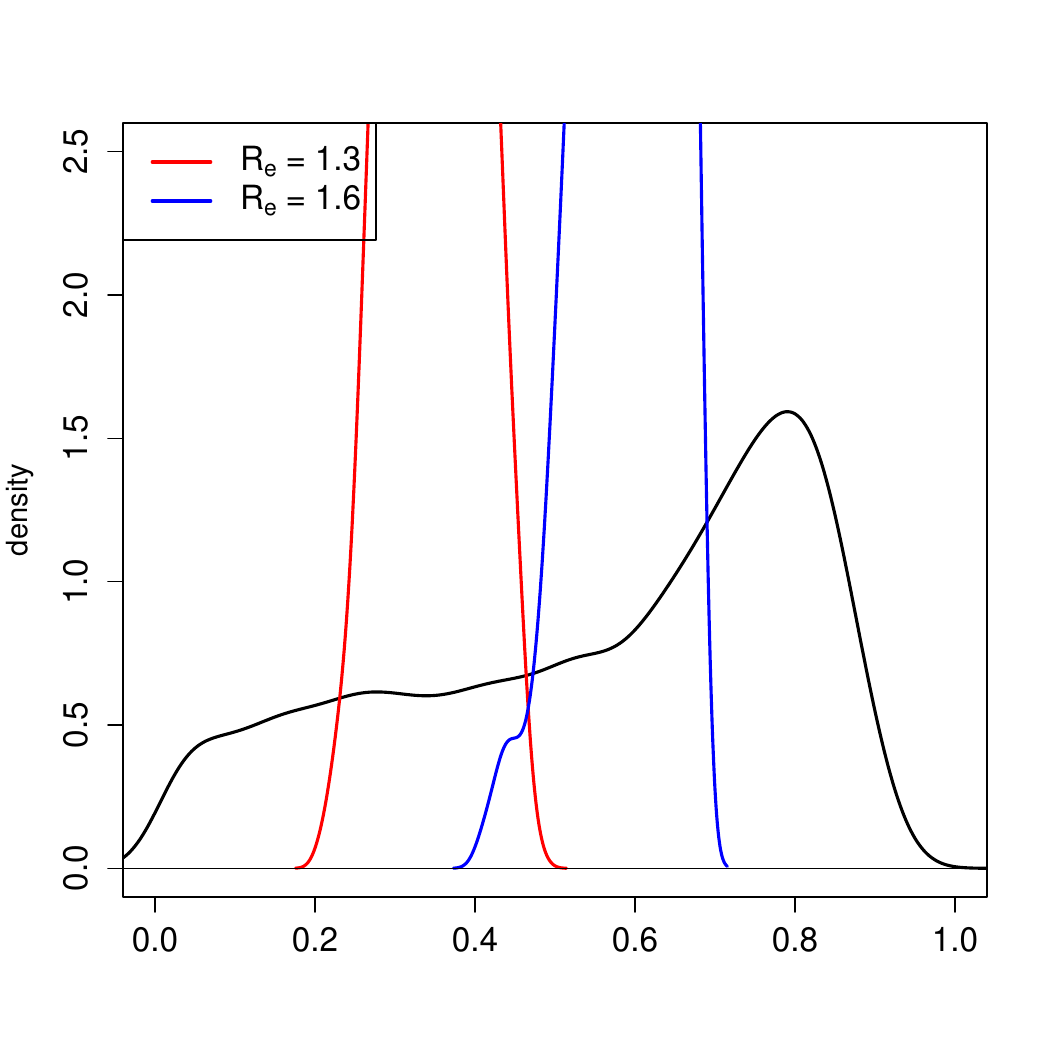}\ 
     \includegraphics[width=0.35 \linewidth]{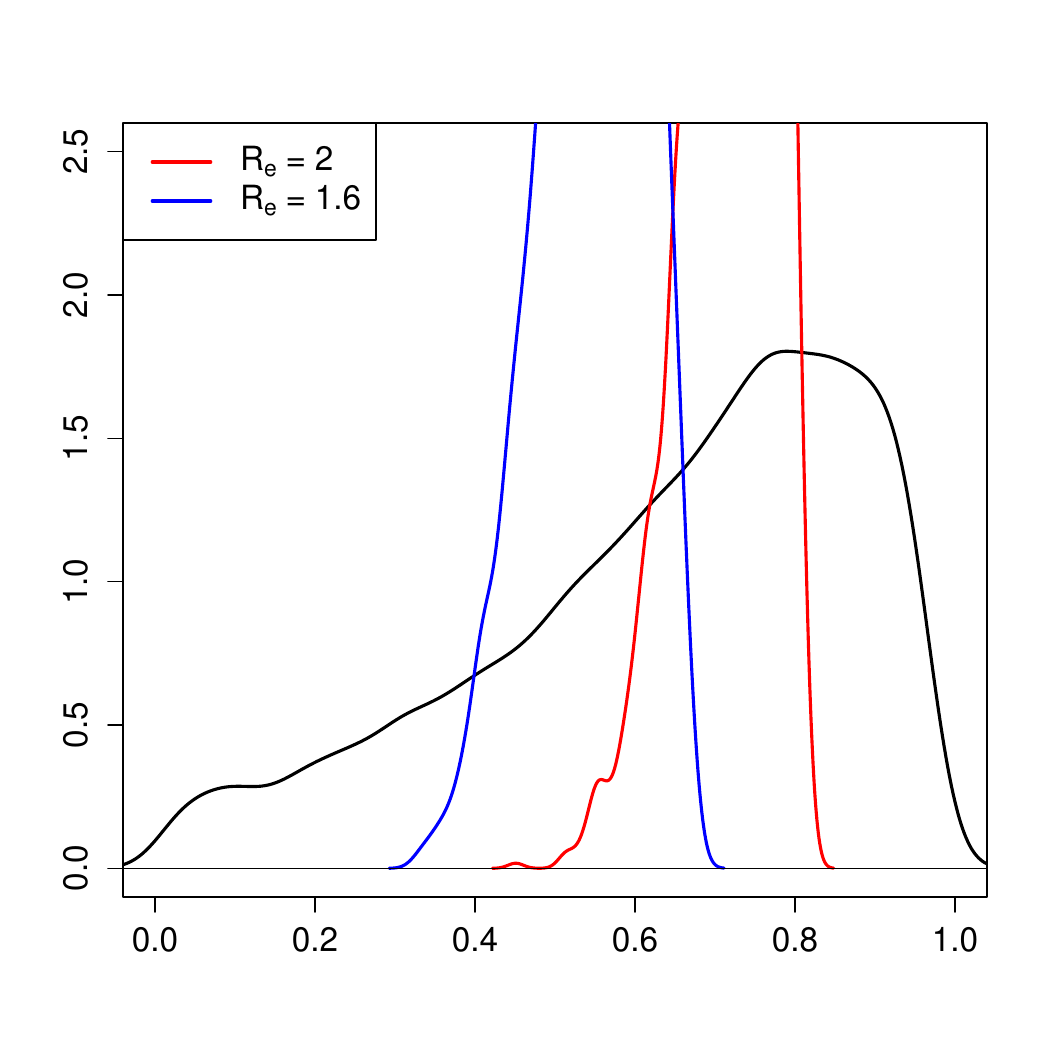}\\[-0.5cm]
       \includegraphics[width=0.35 \linewidth]{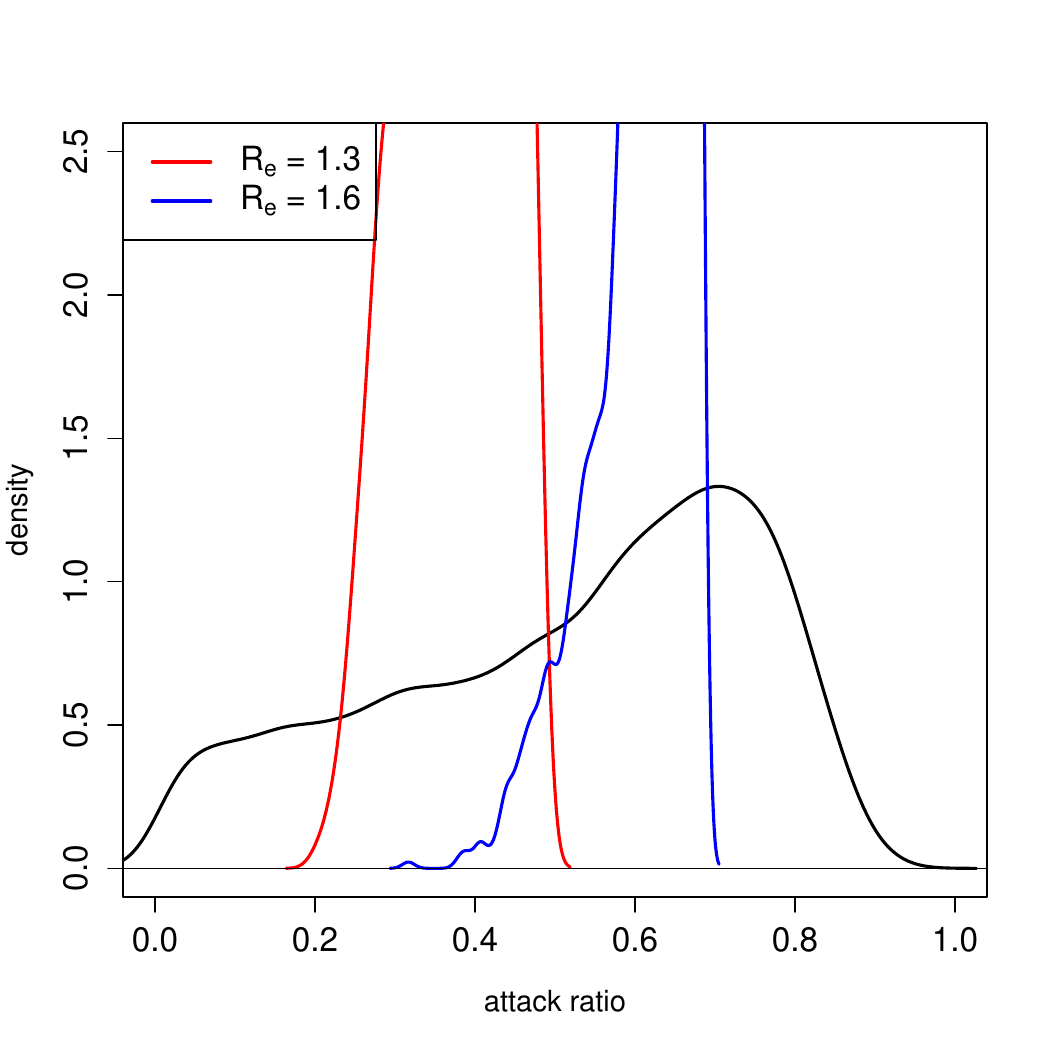}\ 
     \includegraphics[width=0.35 \linewidth]{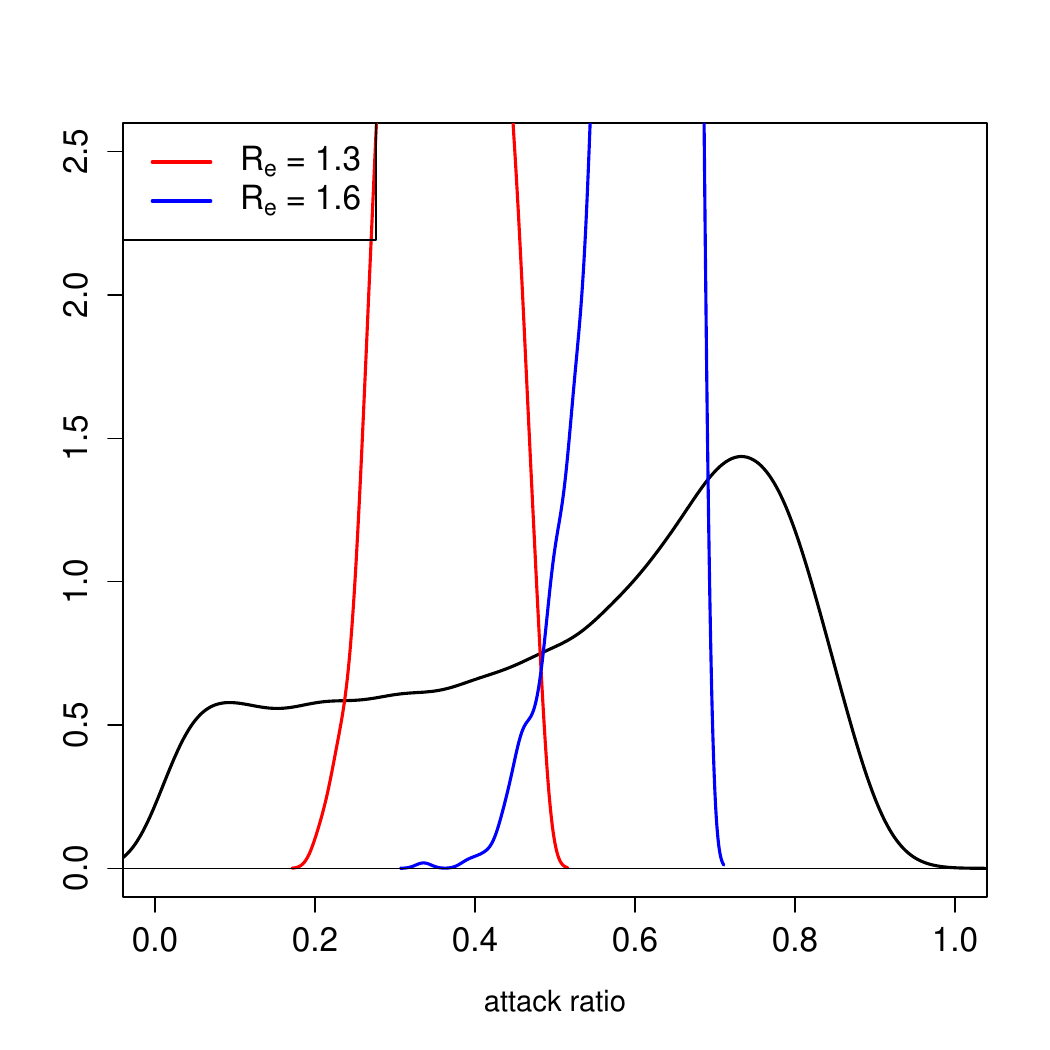}
    \caption{Estimates of the stationary distribution (black curve) of the attack ratio $z^{(k)}$, and of the conditional distributions of $z^{(k)}$ given $R_e^{(k)}=R_e$ (red and blue curves); values used for $R_e$ in the legend.  The four panels correspond to the four distributions of the pairs $(\delta_k,\tau_k)$ as shown in Table \ref{tab:pairs}.}
    \label{fig_statdist}
\end{figure}
In Figure \ref{fig_statdist} we show (black curves) the stationary distribution of $z^{(k)}$ computed through simulation of 20,000 values, smoothed through the \texttt{density} function of R. These are compared to the conditional distributions of $z^{(k)}$ given $R_e^{(k)}$ for  $R_e^{(k)}=1.6$ and $1.3$, being common values for seasonal influenza estimates (e.g.\ \cite{T22}). 

Note that the density shown for the stationary distributions (as well the transition probabilities shown in Figure \ref{fig:dist_givenp}) are defective,
in the sense that the stationary distribution also has probability mass at $z=0$. The value of $\bar \pi_0$ is estimated at 0.25, 0.06, 0.32 and 0.28 in the four cases illustrated in the figure.

Note also that we showed in Subsection \ref{subs:stat} that the conditional distributions of $z^{(k)}$ given $R_e^{(k)}=R_e$ has an atom at $z(R_e)$; this is not so clear in Figure \ref{fig_statdist}, because it is an empirical estimate of the distribution from a simulation through the function \texttt{density}, and because we selected the values of $k$ $R_e^{(k)}$ was close to $R_e$ and not exactly equal to $R_e$ (which is practically impossible). For these reasons, the atom of $z(R_e)$ looks blurred in the Figure.

\begin{figure}[t]
    \centering
    \includegraphics[width=0.35\linewidth]{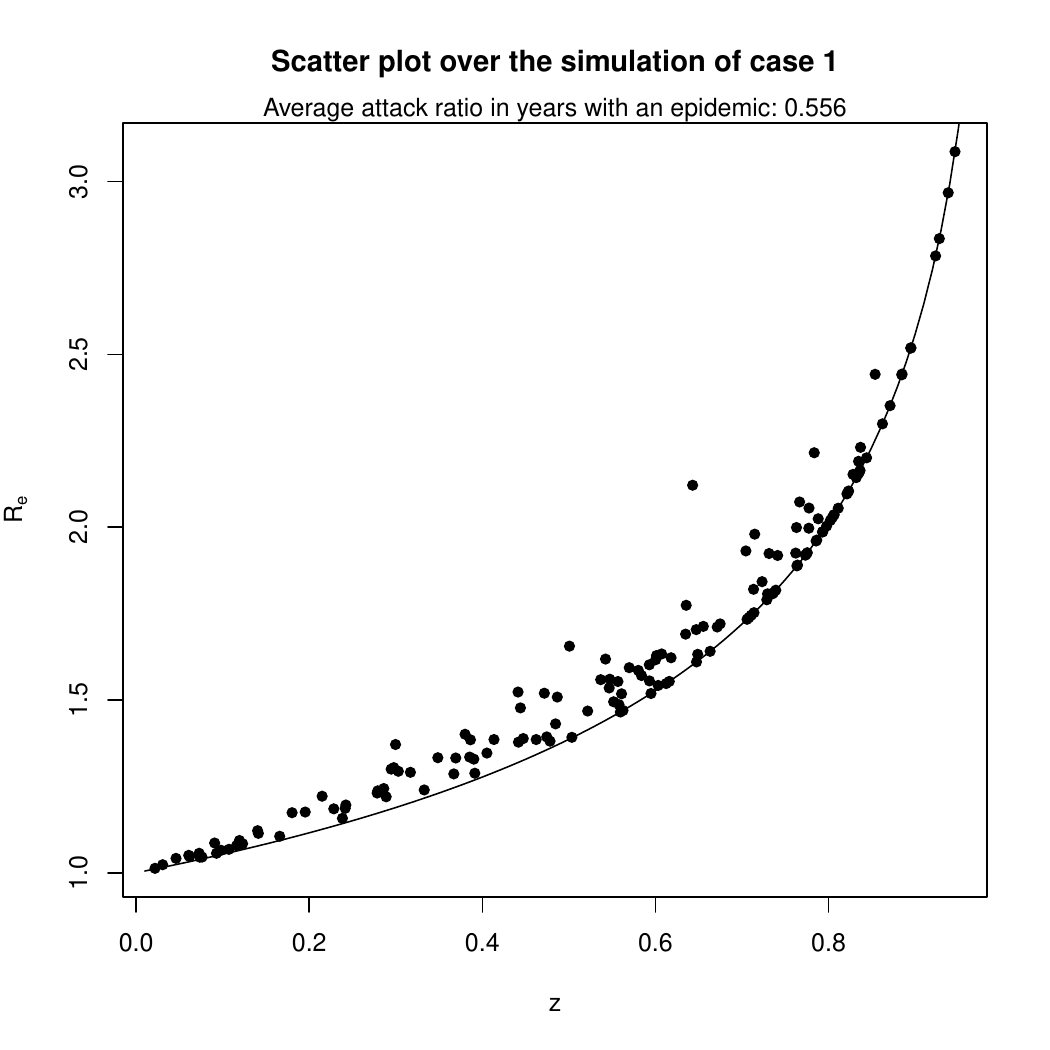}\ 
  \includegraphics[width=0.35\linewidth]{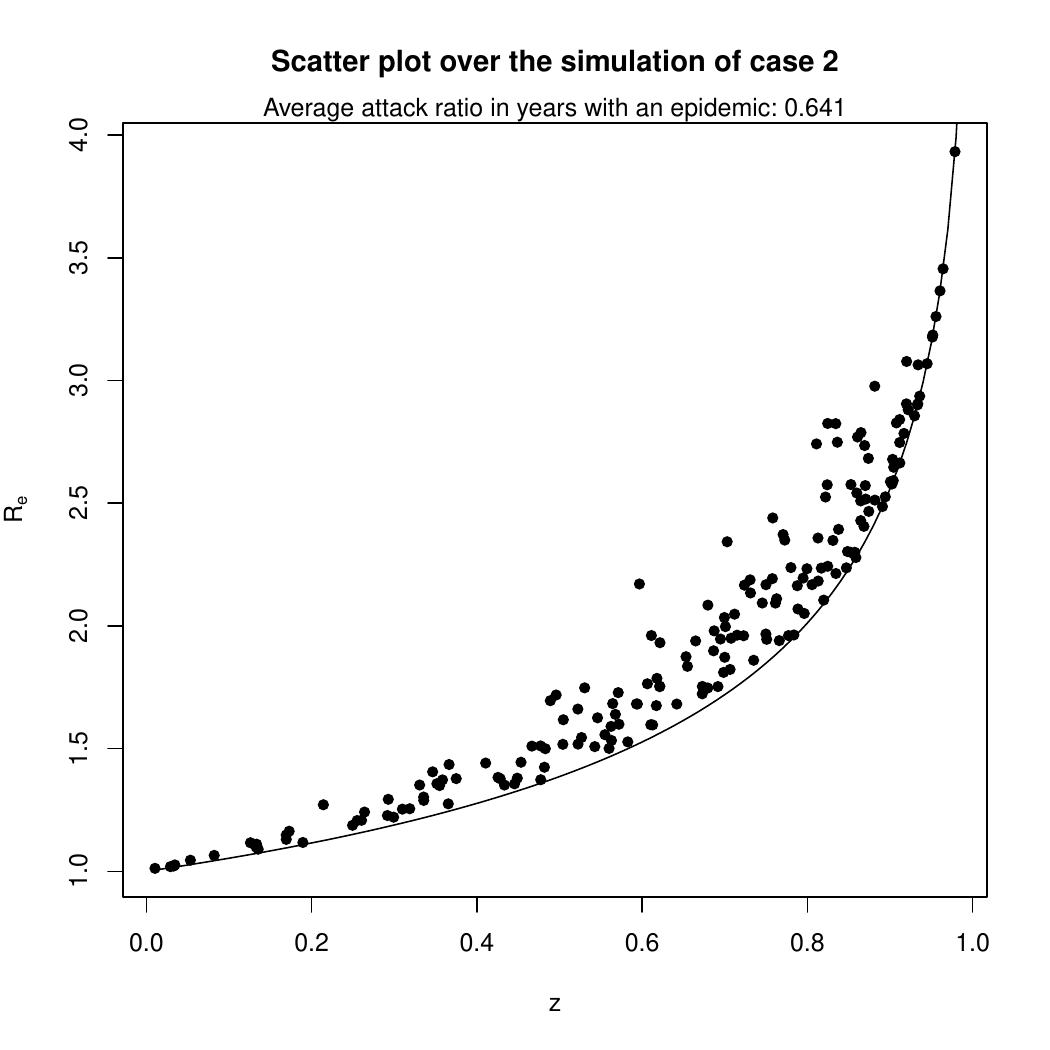}\\
    \includegraphics[width=0.35\linewidth]{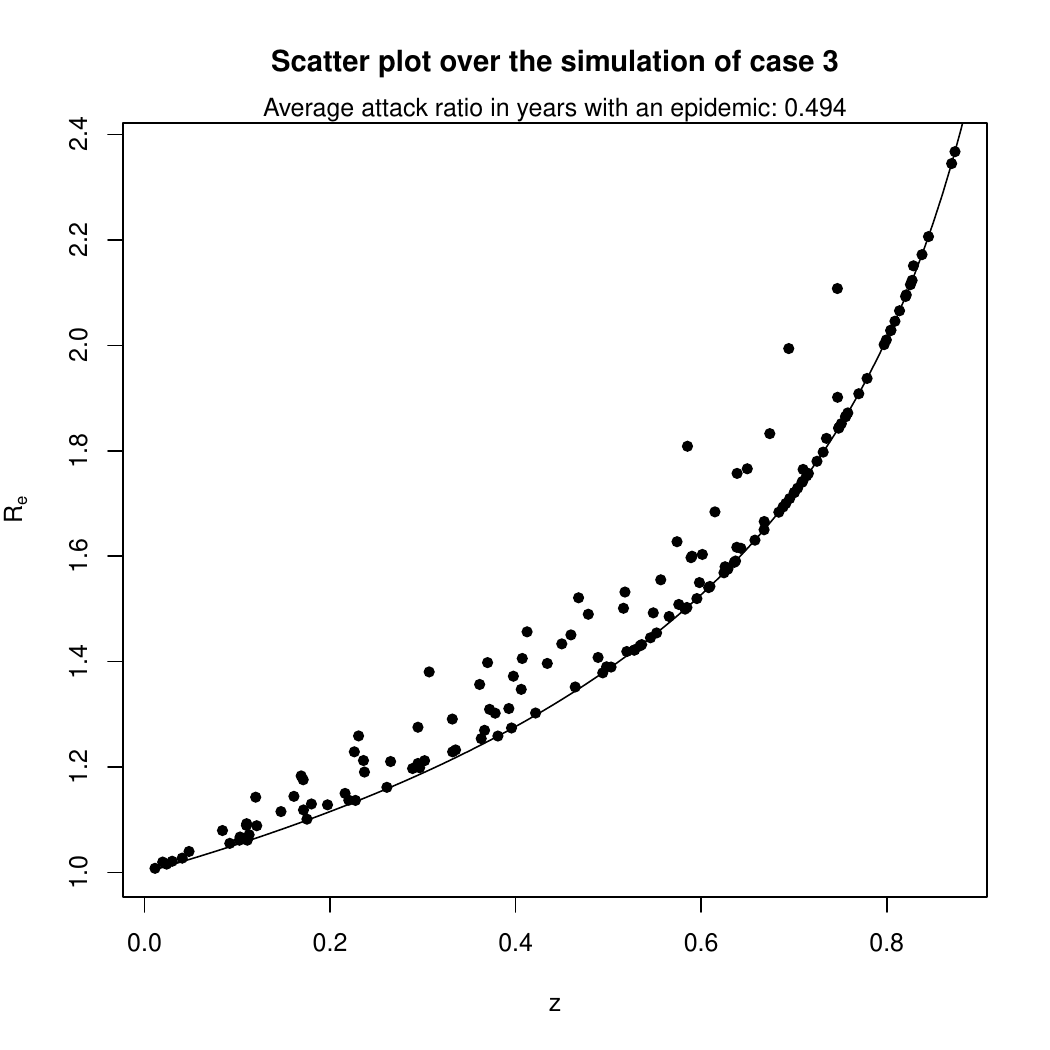}\ 
  \includegraphics[width=0.35\linewidth]{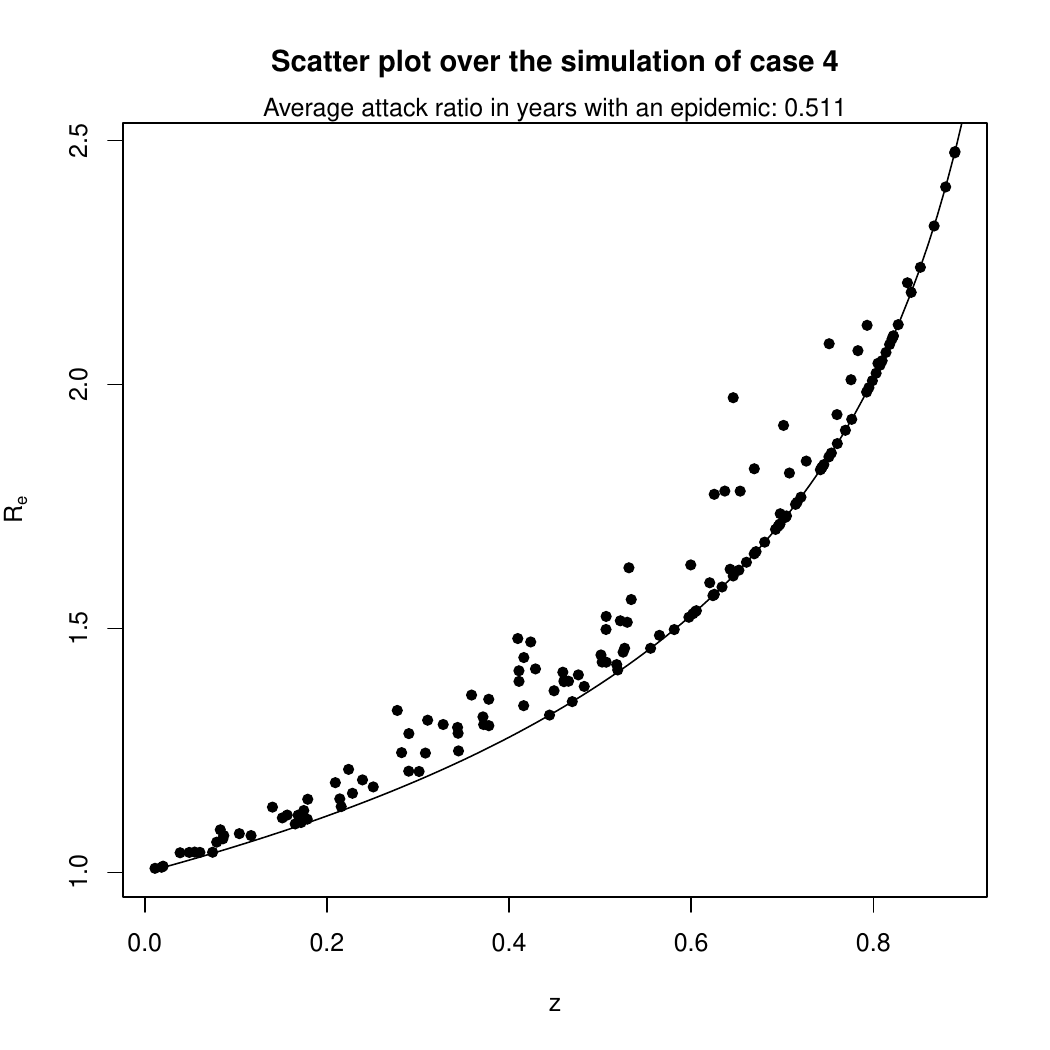}
    \caption{Bivariate graphs of $(R_e^{(k)},z^{(k)})$ found along the simulations of the model with $r=2$. The curves are the functions $R_e = -\log(1-z)/z$ above which the bivariate distributioin always lies, as shown in Section \ref{Sec-r=2}.}
    \label{fig:empcorr}
\end{figure}

In Figure \ref{fig:empcorr} we show the distribution of $(R_e^{(k)},z^{(k)})$ along the simulations. It can be seen that in all cases, the points are rather close to the curve $R_e = -\log(1-z)/z$; in other words,  if the model were correct, predicting the attack ratio from the final size equation using the value of $R_e^{(k)}$ (without using the knowledge of $z^{(k-1)}$) would give a rather good result.

We have also considered the more realistic case of $r=10$, i.e.\ immunity is completely lost 10 years after the last infection; in this case the state space has dimension 17, and thus the transition probabilities are very complicated to compute analytically and visualize. 
\begin{figure}[h]
    \centering
    \includegraphics[width=0.35 \linewidth]{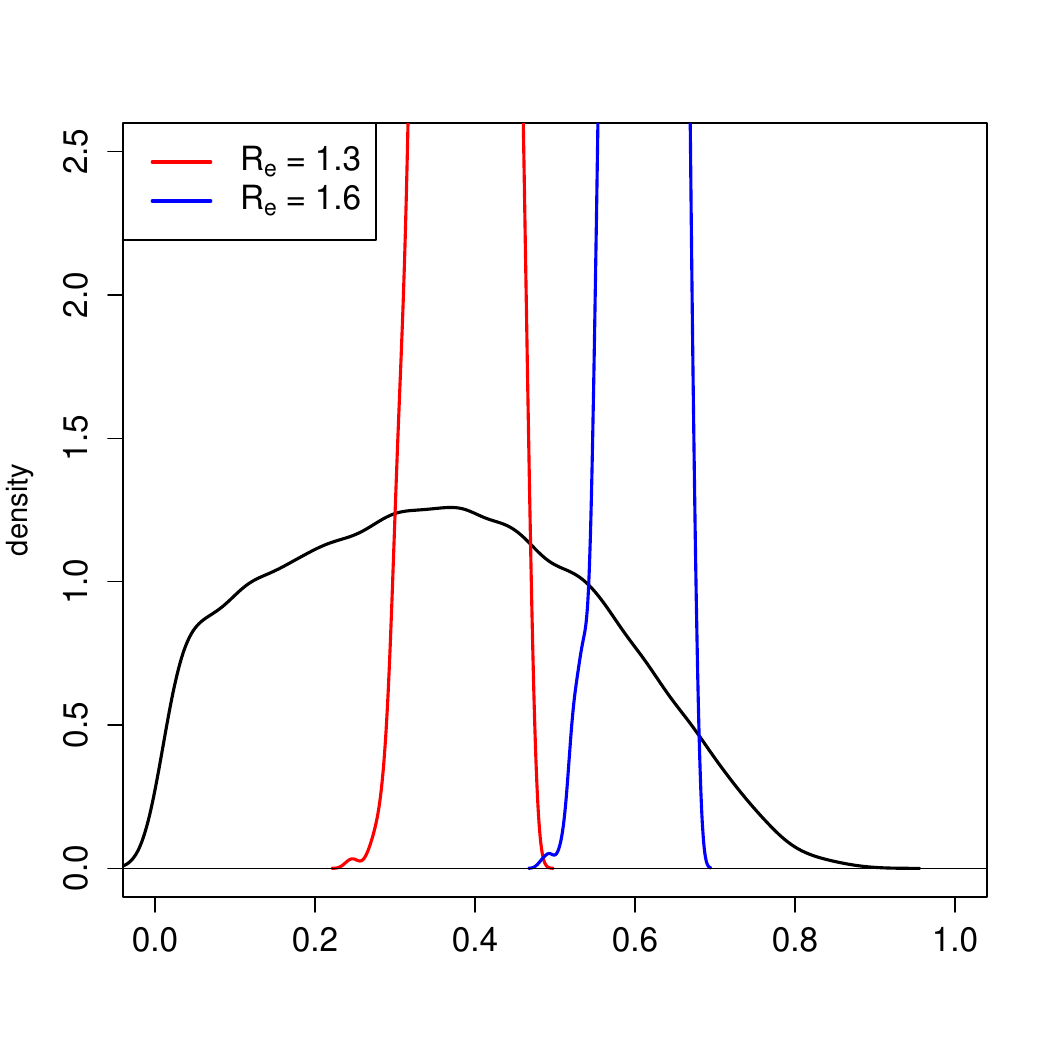}\ 
     \includegraphics[width=0.35 \linewidth]{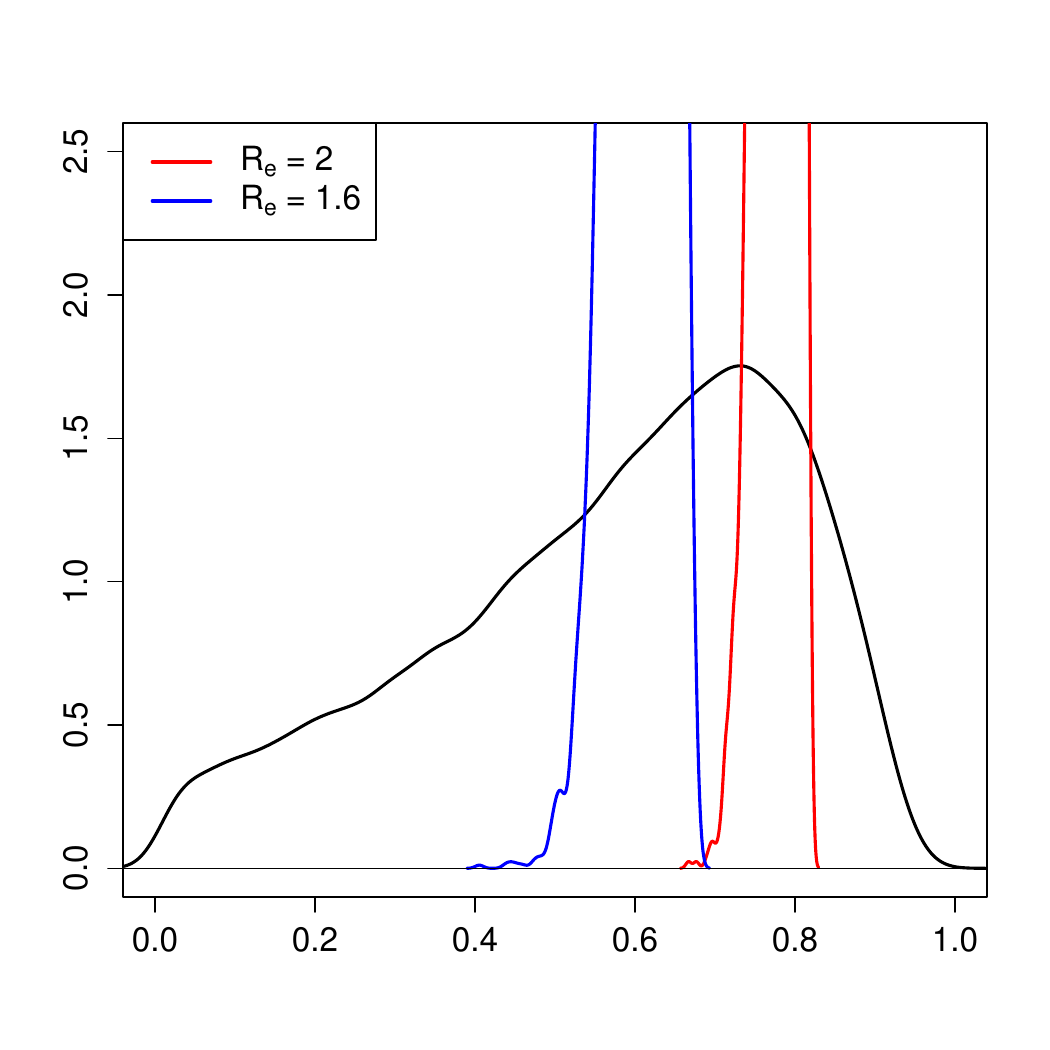}\\[-0.5cm]
       \includegraphics[width=0.35 \linewidth]{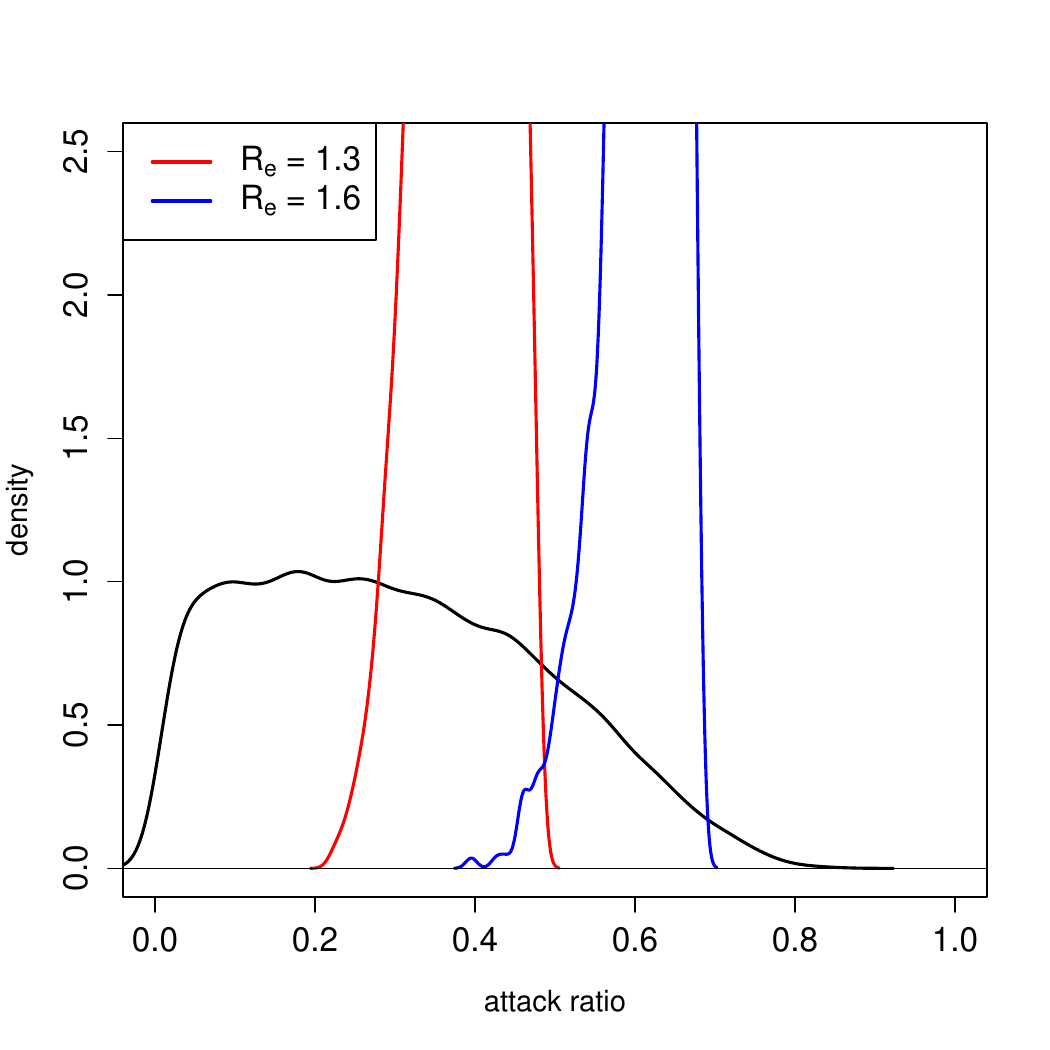}\ 
     \includegraphics[width=0.35 \linewidth]{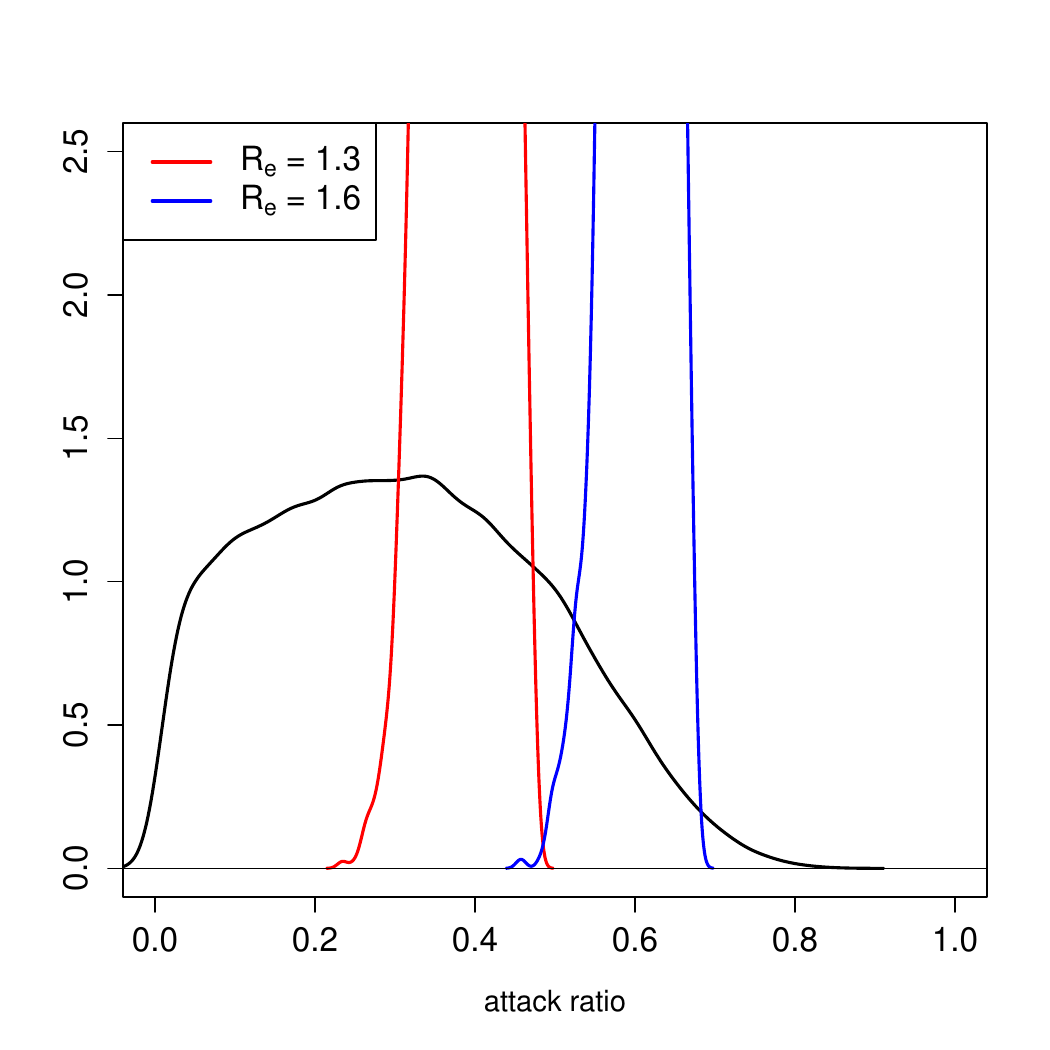}
    \caption{Estimates of the marginal distribution (black curve) of the attack ratio $z^{(k)}$, and of the conditional distributions of $z^{(k)}$ given $R_e^{(k)}=R_e$ (red and blue curves), in the stationary distribution of the model with $r=10$; values used for $R_e$ in the legend.  The four panels correspond to the four distributions of the pairs $(\delta_k,\tau_k)$ as shown in Table \ref{tab:pairs}.}
    \label{fig_statdist_10}
\end{figure}

Therefore, when $r=10$, we limited ourselves to simulations of the model. 
In Figure \ref{fig_statdist_10} we repeat the same procedure used to obtain Figure \ref{fig_statdist}. In this case the stationary distribution is multi-dimensional, and we show estimates of the marginal distribution of the final size $z^{(k)}$. If one compares Figure \ref{fig_statdist_10} to Figure \ref{fig_statdist}, one sees that the densities are shifted to the left with $r=10$, i.e.\ when immunity lasts longer the average attack ratio is lower. On the other hand, once $R_e^{(k)}$ is known, the estimates of $z^{(k)}$ are similar for $r=2$ and $r=10$.

\begin{figure}
    \centering
    \includegraphics[width=0.35\linewidth]{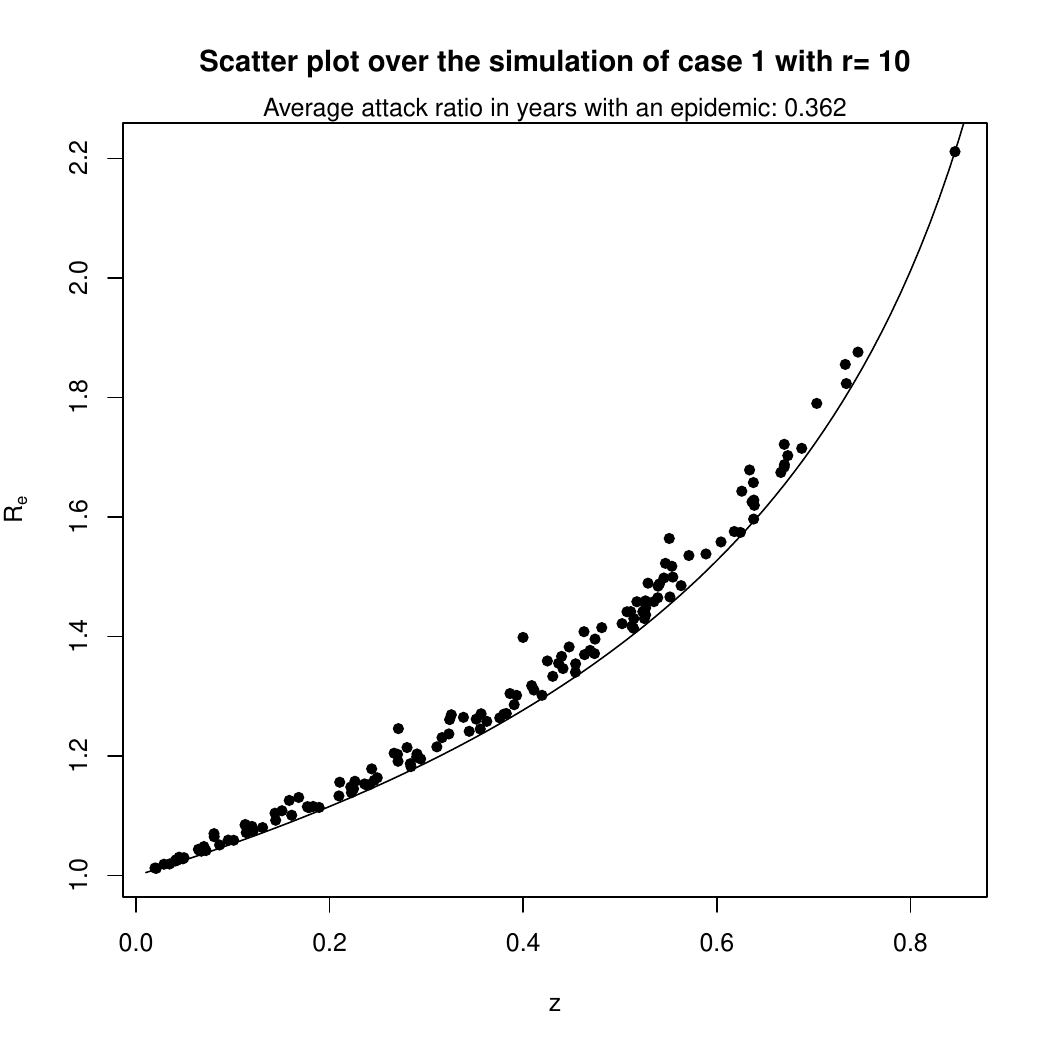}\ 
  \includegraphics[width=0.35\linewidth]{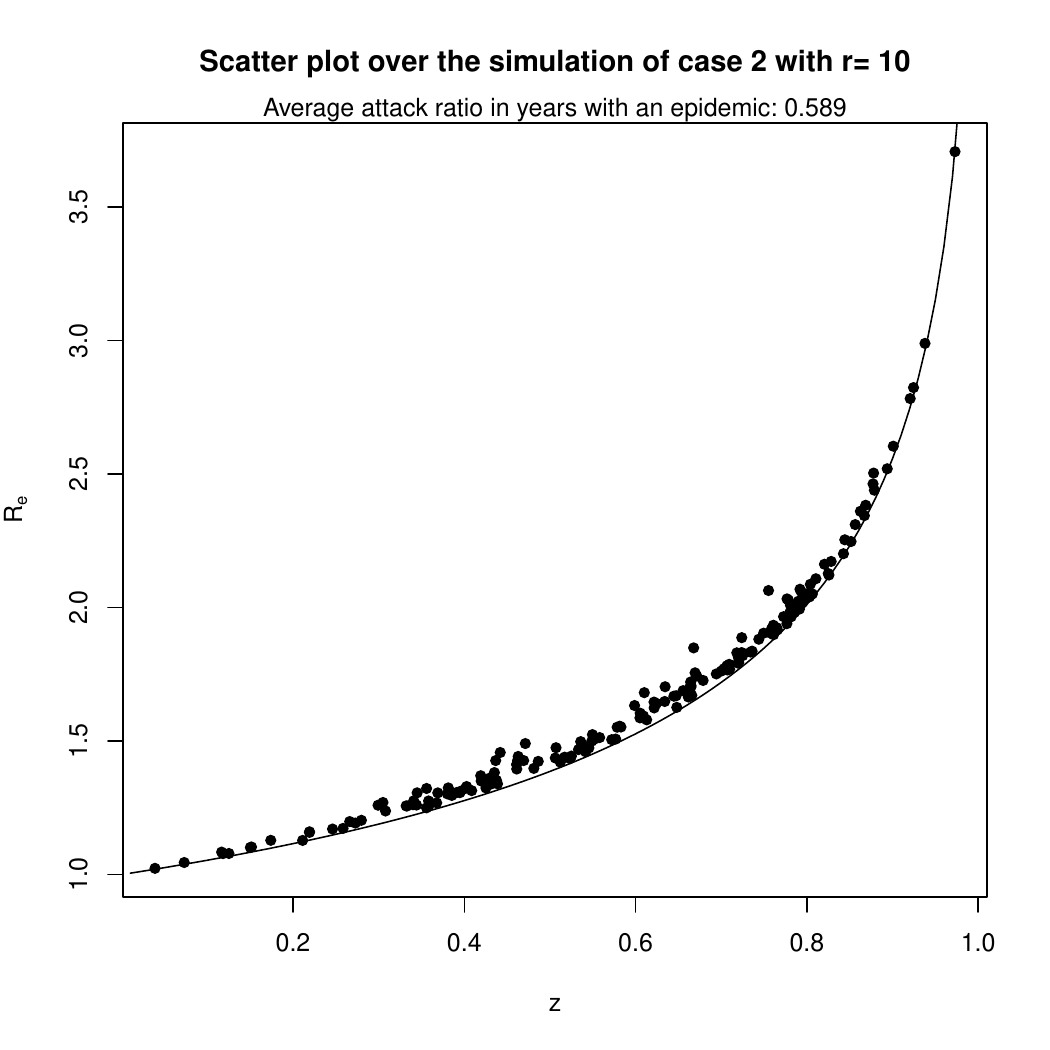}\\
    \includegraphics[width=0.35\linewidth]{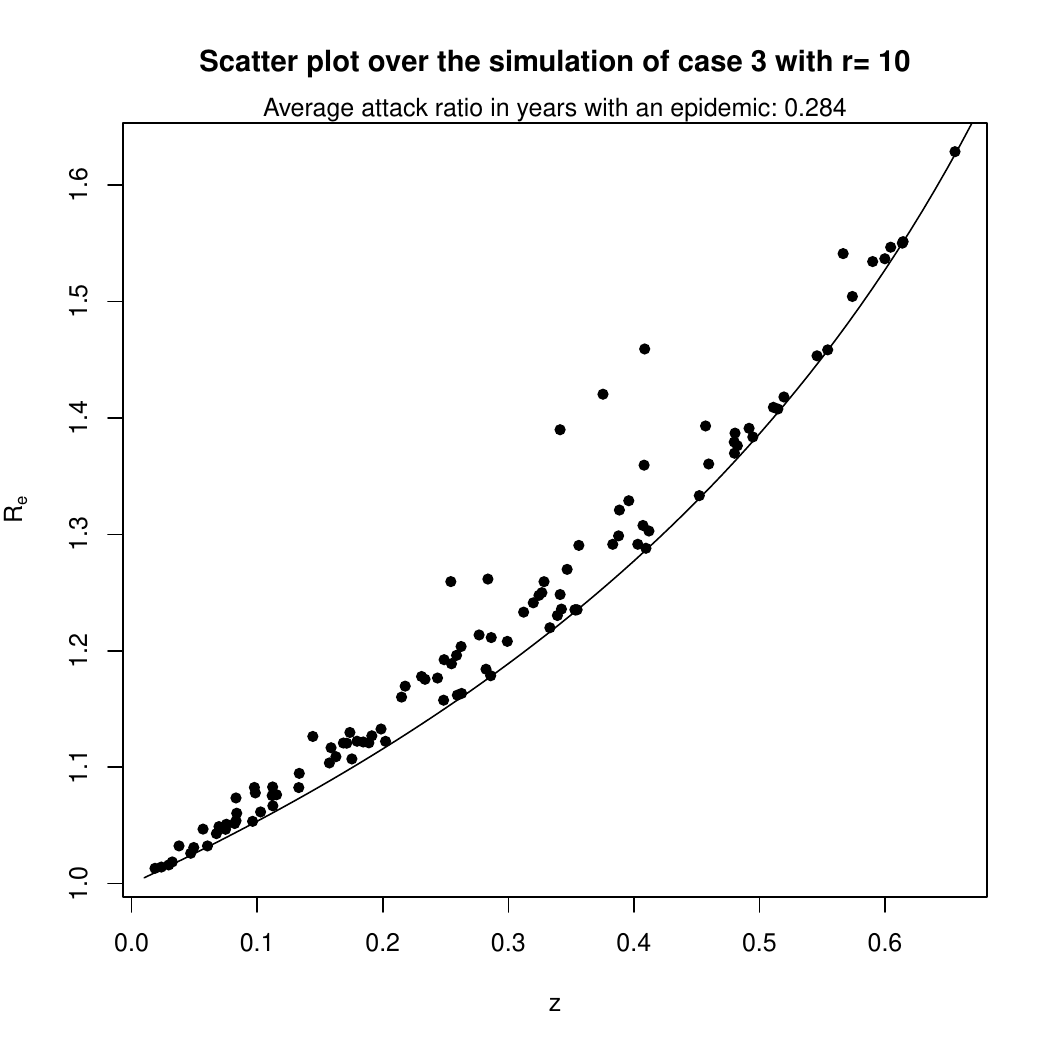}\ 
  \includegraphics[width=0.35\linewidth]{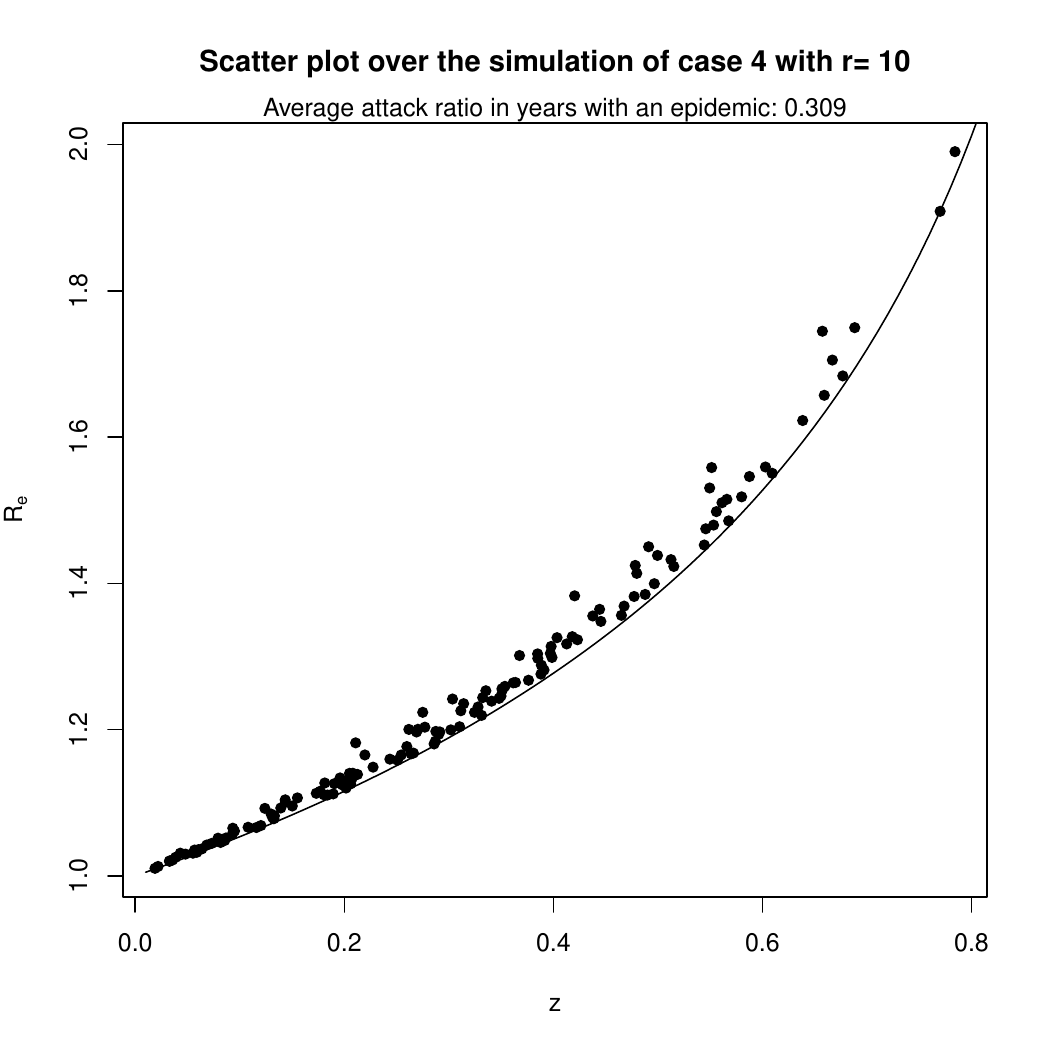}
    \caption{Bivariate graphs of $(R_e^{(k)},z^{(k)})$ found along the simulations of the model with $r=10$. The curves are the functions $R_e = -\log(1-z)/z$ above which the bivariate distributioin always lies, as shown in Section \ref{Sec-r=2}.}
    \label{fig:empcorr10}
\end{figure}

In Figure \ref{fig:empcorr10} we show, analogously to Figure \ref{fig:empcorr}, the pairs $(R_e^{(k)},z^{(k)})$.

Comparing the results obtained in the case $r=10$ with the corresponding ones for $r=2$, one notices that the values of the pairs $(R_e^{(k)},z^{(k)})$ are even closer to the curve $R_e = -\log(1-z)/z$ than for $r=2$. The fraction of years having an epidemic outbreak is similar for $r=10$ and $r=2$, except for case 3, while the average size of the epidemics is lower when $r=10$ than when $r=2$, again less so in case 2, in which immunity drift is generally larger.




\section{Discussion}\label{Sec-Disc}

In the paper we have introduced a stochastic epidemic model for seasonal outbreaks (of e.g.\ influenza), where the virus has a random genetic drift and a new random transmissibility each year. Given the community immunities coming into season $k$, the epidemic outcome, described by the effective reproduction number $R_e^{(k)}$ and final attack ratio $z^{(k)}$, is determined by the drift and transmissibility of this year's virus strain. The novelty of our model lies in considering the yearly drift as random, varying between years, and to also allow for transmissibility to vary between years.

The process of community immunities coming into next year's season, is shown to be a Markov chain having a stationary distribution. The case where immunity only lasts one year is studied in more detail. For this situation, we obtain analytical results for the transition distribution of the Markov chain, its stationary distribution, and also the conditional distribution of the final attack ratio, given the early epidemic growth (corresponding to $R_e^{(k)}$) and incoming immunity. 

The hope is that our model can contribute to the understanding of how and why influenza outbreaks vary in size and initial growth between years, and also potentially help to improve precision in predicting the severity of an influenza season based on observing the start of the outbreak and knowing earlier years outcome. 

Our model of course has several simplifications compared to empirical seasonal epidemics. The perhaps most important extensions for making the model more realistic, are a) to consider more than one circulating strain: influenza currently has three main circulating strains (AH1N1, AH3N2 and B) with separate attack ratios varying between years, and with interactions between the circulasting strains (e.g.\ \cite{Y20}), and b) to allow for heterogeneities in the community, for example in terms of partial immunities and how these affect infectivity. 

In our illustrations we assume that the drift and transmissibility of different years are independent. Our framework however allows for dependecies between years, for example that large genetic drifts rarely happen repeatedly consecutive years, and this can be incorporated, still keeping the Markov property in our model. 

From an applied point of view, it would be interesting to try to fit the model to seasonal outbreaks over a range of years, and to analyse how the fitted stationary distribution, transition distribution, and early season predictions fit to data. It would also be of interest to analyse what underlying distribution of the drift and transmissibility that fits best to observed seasonal outbreak data, which may then be compared to other estimates of yearly drift and transmissibility changes.

Of more mathematical curiosity, it would be interesting to see if a model allowing longer lasting immunity (larger $r$) is more flexible than when immunity only lasts one year.  It would also be interesting to obtain qualitative results regarding concentration of randomness of the final size coming from knowing prior immunity and/or the initial growth rate, in line with what was illustrated numerically in Section \ref{Sec-Numerics}.

\section*{Acknowledgements}
T.B.\ is grateful to the Swedish Research Council (grant 2020-0474) for financial support. This work was initiated when T.B.\ visited University of Trento for 4 weeks in the spring 2025. He is grateful to the hospitality of the Department of Mathematics in Trento. A.P.\ thanks the Italian Ministry for University and Research (MUR) through the PRIN 2020 project ``Integrated Mathematical Approaches to Socio-Epidemiological Dynamics'' (No. 2020JLWP23). A.P.~is member of the ``Gruppo Nazionale per l'Analisi Matematica e le sue Applicazioni" (GNAMPA) of the ``Istituto Nazionale di Alta Matematica" (INdAM).

\section*{Data availability statements}
The paper does not analyse any data. Simulation programmes can be requested from the authors.

\end{document}